\documentclass{elsarticle}
\usepackage[utf8]{inputenc}

\usepackage{subfig}
\usepackage{graphicx}
\usepackage[dvipsnames]{xcolor}
\usepackage{tabularx}
\usepackage{amsmath,amssymb}
\usepackage{mathrsfs} 
\usepackage{gensymb}
\usepackage{cleveref}
\usepackage{mathtools} 
\usepackage{gensymb} 
\usepackage{caption} 
\usepackage{bm} 
\usepackage{cleveref}
\usepackage{numprint} 
\usepackage{multirow} 
\usepackage{multicol} 
\usepackage{longtable}
\usepackage{lineno} 
\usepackage{eufrak}
\usepackage{ulem}
\usepackage{siunitx}
\usepackage{rotating}

\newcommand{\sib}[1]{[\si{#1}]}
\newcommand{\partialder}[2]{\partial_{#2} #1}

\newcommand{\divergence}[0]{\nabla \cdot}

\definecolor{raspberry}{rgb}{1,0,0.2}
\definecolor{darkred}{rgb}{0.5,0.1,0.1}

\newcommand{\m}[1]{\underline#1}
\newcommand{\mm}[1]{\underline{\underline{#1}}}
\newcommand{\qif}[0]{\quad \text{if}}

\graphicspath{{./}} 

\title{A hybrid upwind scheme for two-phase flow in fractured porous media}

\author{Enrico Ballini$^{1,*}$}
\author[1]{Luca Formaggia}
\author[1]{Alessio Fumagalli}
\author{Eirik Keilegavlen$^2$}
\author[1]{Anna Scotti}

\begin{document}

\maketitle
$^1$ MOX, Department of Mathematics, Politecnico di Milano, Piazza Leonardo da Vinci 32, 20133 Milano, Italy. \\
$^2$ Center for Modeling of Coupled Subsurface Dynamics,
Department of Mathematics, University of Bergen, Norway. \\

\noindent
$^*$ Corresponding author, enrico.ballini@polimi.it \\

\noindent
Keywords: finite volume, hybrid upwind, two-phase flow, fractured porous media, mixed-dimensional

\section{Abstract}
Simulating the flow of two fluid phases in porous media is a challenging task, especially when fractures are included in the simulation. Fractures may have highly heterogeneous properties compared to the surrounding rock matrix, significantly affecting fluid flow, and at the same time hydraulic aperture that are much smaller than any other characteristic sizes in the domain. Generally, flow simulators face difficulties with counter-current flow, generated by gravity and pressure gradients, which hinders the convergence of non-linear solvers (Newton).

In this work, we model the fracture geometry with a mixed-dimensional discrete fracture network, thus lightening the computational burden associated to an equi-dimensional representation. We address the issue of counter-current flows with appropriate spatial discretization of the advective fluid fluxes, with the aim of improving the convergence speed of the non-linear solver. In particular, the extension of the hybrid upwinding to the mixed-dimensional framework, with the use of a phase potential upstreaming at the interfaces of subdomains.

We test the method across several cases with different flow regimes and fracture network geometry. Results show robustness of the chosen discretization and a consistent improvements, in terms of Newton iterations, compared to use the phase potential upstreaming everywhere.


\section{Introduction}
\label{sec:introduction}
Fluid flow and multiphase transport in fractured porous media is of critical importance for subsurface engineering: fractures may form major pathways for fluid flow,
indeed, as fractures can have significantly higher permeability than the surrounding host rock,
they may constitute the main pathways for fluid flow.
Rapid flow through fracture networks can be a desired effect, for instance in the
production of geothermal energy from hard rocks, but it may be detrimental for the
storage of carbon dioxide and nuclear waste \cite{iding2010evaluating,tsang2015hydrologic}.
Moreover, for energy storage and production, fluid exchange between fractures and the host rock is important,
while transport in fracture network also play an important part in natural subsurface flows,
including thermal convection, e.g., \cite{lister1974penetration,fujiwara2003crustal}.

Numerical simulations are valuable tools to study and understand flow in fractured porous media,
however, construction of adequate simulation models is challenging. 
The contrast in permeability between fractures and the host rock,
together with the lack of scale separation in fracture length \cite{bonnet2001scaling}, imply that traditional upscaling is difficult for fluid flow,
and even more so for transport processes.
To classify the many models that have been developed to meet this challenge, it is useful to
consider whether the fractures are represented explicitly or by equivalent continua \cite{Berre2019b}.
In equivalent continuum models, the fractures and host rock are represented by one or multiple overlapping domains,
with fluid flow taking place both within and possibly between the domains.
Such continuum models can preserve the heterogeneity in flow properties between fractures and the host rock,
however, estimating the transfer coefficients between the media may be challenging \cite{karimi2006generation,gong2008upscaling,geiger2013novel}.
Nevertheless, compared to alternative methods, continuum models have a relatively low computational cost, and variants thereof
have been applied to large-scale simulations.
In so-called embedded discrete fracture matrix models (EDFM), fractures are explicitly represented at the continuous level, but not in the computational grid,
see for instance \cite{lee2001hierarchical,fumagalli2016upscaling,ctene2017projection}. 
This results in simulation models that, in a sense, are very similar to that of continuum models, 
however, the explicit representation of fractures 
eases the calculation of the flow exchange between fractures and host rock.
EDFM models have been developed to a high level of sophistication, see for instance \cite{li2008efficient}.

Our main interest herein is in so-called discrete fracture matrix (DFM) models, which represent fractures explicitly in the computational grid.
DFM models usually represent the fractures as lower-dimensional objects embedded in the host domain,
resulting in a mixed-dimensional geometry \cite{Martin2005,Mishra2010,Jaffre2011,Fumagalli2012d,Ahmed2015,Brenner2015,Brenner2016a,Ahmed2017,Brenner2018,Boon2018}.
In these models, equations and constitutive laws can be represented in the host rock, in fractures, and on the rock-fracture interface.
The explicit representation of fractures can complicate grid construction and require a high number of grid cells, 
and thereby limit the domain size and number of fractures that can be included.
Nevertheless, the detailed representation of the fracture network geometry makes DFM models ideally suited
to study the interaction between physical processes such as flow and transport \cite{Fumagalli2012a,Fumagalli2020f,Berre2021,berge2024numerical},
and also mechanical deformation of fractures \cite{garipov2016discrete,stefansson2021fully,novikov2022scalable}.
Particularly relevant to this work is the discretization of multiphase flow, 
which has been reported in, for instance, \cite{Fumagalli2012d,Brenner2015,glaser2017discrete,Brenner2018}.


Multiphase flow in porous media can be described by a pressure equation which, in the assumption of negligible capillary pressure,
has an elliptic character, combined with a set of transport equations that are essentially hyperbolic in nature \cite{Trangenstein1989}.
For practical simulations, the equations are usually discretized by finite volume methods \cite{rasmussen2021open}.
In particular, to avoid unphysical oscillation at the discrete level, phase mobilities have traditionally been discretized by upwinding each phase individually,
using the so called phase potential upstreaming (PPU) \cite{brenier1991upstream}.
PPU is monotone and first order convergent, and often produces results that are reasonably accurate. 
However, for flow governed by a mixture of viscous and gravitational forces, the upstream direction 
assigned to each individual grid face in the discrete model 
is prone to flipping during Newton iterations, leading to convergence issues for the non-linear solver, e.g. \cite{li2015nonlinear}.

To overcome these convergence problems, several improvements of nonlinear solvers can be used, 
including trust-region \cite{wang2013trust,moyner2017nonlinear} and reordering methods \cite{kwok2007potential,natvig2008fast}.
Of interest to us herein is an approach that replaces the PPU treatment of mobilities with a discretization which is in a sense smoother,
and thus less prone to changes in the upstream direction. 
Among the possible techniques to achieve such regulatization, 
we focus our attention on the method known as hybrid upwinding (HU), which was developed in a series of papers \cite{Lee2015,Hamon2016,Hamon2016a,alali2021finite,Bosma2022}.
At the core of the HU approach is the representation of the viscous flow by a total velocity field,
with transport of individual phases taken as gravitational deviations from the total velocity.
The total velocity is discretized with weighted averages of mobilities. 
Compared to PPU, HU posseses enhanced smoothness which leads to significant improvements
of the performance of the Newton solver, to the price of somewhat increased numerical diffusion.

In this work, we consider two phase flow for fractured porous media, and extend the HU 
approach to mixed-dimensional DFM models with explicit representation of fractures in the 
computational grid.
Noting that HU applied to fracture networks in a capillary dominated regime has already been reported \cite{alali2021finite},
we limit ourselves to the case of negligible capillary forces, and focus on the impact of
viscous and gravitational effects.
Compared to standard porous media, the presence of fracture networks introduces some additional difficulties,
mainly the coupling among domains of different dimensions and the strong contrast in permeability, 
resulting in faster dynamics in the fracture network.
Through a series of numerical experiments in two- and three-dimensional domains, we show
that the extension of HU consistently outperforms PPU in terms of performance of a Newton solver.
The computational gains increase with simulation complexity, indicating that our extension
can be a key ingredient in enabling two-phase 3D DFM simulations in regimes where
gravitational effects play an important part.

The paper is structured as follows: The governing equations are presented in Section \ref{sec:mathematical_model}, 
while in Section \ref{sec:discretization}, we present the discretization methods with emphasis on the different approaches to upstreaming. 
Numerical tests are presented in Section \ref{sec:numerical_validation}, while Section \ref{sec:conclusion} contains concluding remarks.

\section{Domain representation} \label{sec:domain}

Due to the complexity of the mixed-dimensional representation, we 
devote this section to introduce and describe such framework. Starting from \Cref{subsec:domain1}, we introduce some
basic notations when the domain does not contain fractures. We extend the concepts in \Cref{subsec:domain2} when fractures are presents and in \Cref{subsec:domain3} we discuss the strategy adopted for the coupling between objects of different dimensions.

\subsection{Fracture-less domain}\label{subsec:domain1}

Let  $\Omega\subset \mathbb{R}^d$, with $d=2$ or $3$, to be a Lipschitz continuous domain representing the porous media when fractures are not present, with boundary $\partial \Omega$ and outer unit normal
$\upsilon_{out}$. Being the problem time dependent, we consider the time interval $(0, T]$, with $T>0$ the final time in \sib{\second}, and introduce the space-time domain $\Omega^T = \Omega \times (0, T]$.

\subsection{Fractured domain}\label{subsec:domain2}

We consider here the approach described in \cite{Berre2019b}, where fractures are explicitly represented. Even for larger fractures, their thickness or aperture is orders of magnitude smaller than their typical lateral extensions.
We thus approximate fractures with lower dimensional objects immersed in the rock domain. For simplicity, we assume that fractures are planar objects and we indicate with $\varepsilon$ their thickness, in \sib{\meter}, assumed to be constant for each fracture.

\begin{figure}
    \centering
    \subfloat[2D and 1D domain.]{\includegraphics[width=0.5\textwidth]{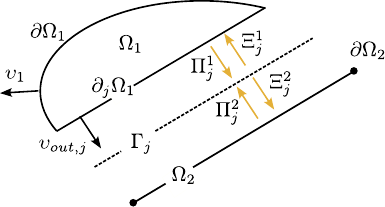}}
    \hspace{0.5cm}
    \subfloat[Y-shape intersection of 1D domains.]{\includegraphics[width=0.35\textwidth]{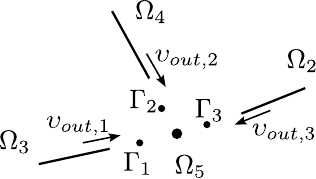}}
    \caption{Mixed-dimensional domain. (a) Domains $\Omega_1$ and $\Omega_2$ are connected through the mortar interface $\Gamma_j$. (b) Example of three branches, in 3D fracture intersections and in 2D fractures, that intersect to a point $\Omega_4$.}
    \label{fig:mixed_dim_domain}
\end{figure}

The modelling of the fractures leads to subdividing the whole domain $\Omega\subset\mathbb{R}^d$ into subdomains of different physical dimensions $\Omega_i$, $i=1,\ldots,I$, with $I$ the total number of subdomains, such that $\overline{\Omega} = \cup_i \overline{\Omega}_i$ and $\Omega_i\cap{\Omega}_j=\emptyset$ for $i\ne j$. Each domain might represent the porous media, a fracture or a fracture intersection. For example, if $d=3$, the subdomains representing the rock matrix have dimension $3$, each reduced fracture  has dimension $2$, fracture intersections have dimension $1$, and the possible intersection of intersections dimension $0$. If $d=2$, then all the aforementioned dimensions should be scaled by $1$. See Fig.\ref{fig:mixed_dim_domain} for an example of a mixed-dimensional domain.

\subsection{Interface coupling}\label{subsec:domain3}
Flow in the subdomains $\Omega_i$ is interconnected through flux exchanges across the interface between domains with codimension equal to 1. We remark that there is no direct interaction between subdomains with codimension greater than 1, such as a 3D matrix rock and 1D intersection of fractures.  We call these fluid fluxes \textit{mortar fluxes} and can be interpreted as Lagrange multipliers that enforce the correct mass balance between subdomains \cite{Boon2018,Nordbotten2018}. In our formulation, we indicate the interface between domains, which we call \textit{mortar interfaces} (or simply mortars) explicitly as $\Gamma_j$, $j = 1, \ldots, J$, with $J$ the total number of mortars. This will be useful for the discretization approach presented in \Cref{sec:discretization}.
See Fig.~\ref{fig:mixed_dim_domain}(a) for an illustration of a 2D matrix domain with one line fracture and a mortar domain interfacing the former domains. 

On each subdomains $\Omega_i$ and mortar $\Gamma_j$ the variables and data are marked by the subscript $i$ or $j$, respectively.
We indicate the mapping from the  $j$-th mortar $\Gamma_j$ to the boundary of the related $i$-th subdomain $\Omega_i$ by $\Xi_j^i$.
We also introduce a map from a $j$-th mortar $\Gamma_j$ to a neighbouring $\Omega_i$, denoted by the symbol $\Pi_j^i$, see Fig.~\ref{fig:mixed_dim_domain}. 
These maps are relevant for the discrete problem, in particular when the meshes are non-conforming across subdomains and mortars, see \Cref{sec:discretization}.

For a subdomain $\Omega_i$, the set of neighbouring mortars is split into mortars, denoted by $\hat{S}_i$, that connect $\Omega_i$ to subdomains of higher dimension, and mortars that connect $\Omega_i$ to subdomains of lower dimensions, represented by the set $\check{S}_i$.
Conversely, we denote by $\hat{R}_j$ the set of subdomains facing mortar $\Gamma_j$ with $\dim(\Omega_i) > \dim(\Gamma_j)$, and with $\check{R}_j$ the set of subdomains facing the mortar $\Gamma_j$ with $\dim(\Omega_i) < \dim(\Gamma_j)$. In the following, we use the abbreviation $l$ to denote the indices of domains $\Omega_i \in \check{R}_j$, and, analogously, $h$ to denote the indices of domains $\Omega_i \in \hat{R}_j$. Thus, for example, $\Pi_j^h$ and $\Pi_j^l$ are the maps from,  respectively, $\Omega_i \in \hat{R}_j$ and $\Omega_i \in \check{R}_j$ towards the mortar $\Gamma_j$.

The normal at the boundaries is $\upsilon_i$, whereas $\upsilon_{out, j}$ is the unit vector at the boundary of the higher dimensional domain, pointing outwards towards $\Gamma_j$.
$\mathcal{N}_j$ is a map from $\overline{\Omega}$ to $\Gamma_j$ such that, given $u : \overline{\Omega} \rightarrow \mathbb{R}$, we have $\mathcal{N}_j(u) : \overline{\Omega} \rightarrow \Gamma_j$. Its relevance will become apparent when we deal with the discretization procedure, where we give a more precise definition.
With $\partial_j \Omega_i$ we indicate the boundary of $\Omega_i$ in contact with the mortar domain $\Gamma_j$. 
We denote the subdomain codimension extension with $\varepsilon_i^{a_i}$, where $a_i = d - \dim(\Omega_i)$. 
For example, the codimension extension (thickness) of a planar fracture immersed in a 3D rock, would be $\varepsilon^{3-2} = \varepsilon$, while the codimension extension (area) of the intersection of two planar fractures would be $\varepsilon^{3-1} = \varepsilon^2$. 
Similarly, we define $b_j = d-\dim(\Gamma_j)$. 

\section{Mathematical model}
\label{sec:mathematical_model}

In this part we introduce the mathematical model considered for the 
two-phase flow. First, in \Cref{sec:equi-dimensional} we present the model for a continuous medium without fractures. Subsequently, in 
\Cref{subsec:dimensional} we discuss a dimensional analysis and obtain dimensionless groups, useful to set up the simulations in the examples. Finally, in \Cref{sec:mixed-dimensional} we present the mixed-dimensional model, an approach to approximate the fractures in according to their features.

\subsection{Continuous model}\label{sec:equi-dimensional}

We assume no fractures and that in the void spaces of the porous media coexists two phases that fulfil the following assumptions: they are immiscible, isothermal, and non-reactive, with null capillary pressure. 
For their description, we consider the classical two-phase flow model in porous media, see \cite{Bear1972, Chavent1986, Helmig1997, Nordbotten2011} for more details. In the sequel, we will indicate with a subscript $0$ and $1$ data and variables associated to each phase.

The primary variables we are considering in our model are the saturations $S_\ell :\Omega^T \rightarrow [0, 1]$, that are dimensionless, for each phase $\ell=0,1$, and the pressure $p:\Omega^T \rightarrow \mathbb{R}$, in \sib{\pascal}, equal for the two phases because of the assumption of null capillary pressure.
The porous medium is characterized by the following properties: the porosity $\phi$, that is dimensionless, and the intrinsic permeability of the rock $K$, in \sib{\square\meter}. Each fluid phase is characterized by:
 the density $\rho_\ell$, in \sib{\kilogram\per\cubic\meter},  the dimensionless relative permeability $k_{r, \ell} : [0, 1]\rightarrow [0,1]$, dependent on the phase saturation $S_\ell$, the dynamic viscosity $\mu_\ell$, assumed to be constant, in \sib{\pascal\cdot\second}. To simplify the notation, we 
introduce the phase mobility $\lambda_\ell = k_{r, \ell} / \mu_\ell$, in \sib{\per\pascal\per\second}.
We set $g$ to be the gravity field constant assumed to be equal to $9.81 \sib{\meter\per\square\second}$, and we consider the vertical coordinate $z$ pointing upwards, so the gravity vector is $-g\nabla z$. Finally,
$f_\ell$ is a source or sink term associated to each phase $\ell$, in \sib{\kilogram\per\cubic\meter\per\second}. 

The problem is to find $(S_0, S_1, p)$ such that for $\ell=0,1$ we have
\begin{subequations}\label{eq:mass_bal_equi}
\begin{align}\label{eq:mass_bal_equi_eq}
    \begin{aligned}
        &\partialder{\phi(\rho_\ell S_\ell)}{t} + \divergence Q_\ell = f_\ell\\
        &Q_\ell = -\rho_\ell \lambda_{i} K \left(\nabla p + \rho_\ell g \nabla z \right)
    \end{aligned}
     & \quad \text{in } \Omega^T,
\end{align}
where $Q_\ell$, in \sib{\kilogram \per \square \meter \per \second}, is the mass flux of each phase, 
which is proportional to the gradient of the phase potential, 
\begin{equation}\label{eq:phase_potential}
    \Phi_\ell = p + \rho_\ell gz.    
\end{equation}
Associated to the previous equations, in our numerical experiments we consider the following initial and boundary conditions
\begin{align}\label{eq:bc_ic}
    \begin{aligned}
    &S_\ell(t = 0, x) = \overline{S}_\ell(x) && \quad \text{in } \Omega,\\
    &p(t=0, x) = \overline{p}(x)&& \quad \text{in } \Omega,\\
    &Q_\ell\cdot \upsilon_{out} =0 
    && \quad \text{on } \partial \Omega \times (0, T],
    \end{aligned}
\end{align}
where $\overline{S}_\ell \in [0,1]$, in \sib{\cdot}, and $\overline{p}$, in \sib{\pascal}, are given functions representing the initial values for the two phase saturation and pressure, respectively. We assume the following constraint for the saturations 
\begin{gather}\label{eq:saturation_eq}
    S_0 + S_1 = 1 
    \quad \text{in } \Omega^T.
\end{gather}
We remark that we need to obey compatibility conditions between the data defining the problem. In particular the initial value of the saturation has to respect the contraint $\overline{S}_0 + \overline{S}_1 = 1$, and, for the given boundary conditions, steady state can only be reached when the source terms, $f_\ell$, have zero average.

To close the system, we consider a constitutive equation to relate the phase 
density, assumed to be a liquid and thus nearly incompressible, with the pressure, and a model for the relative permeability, namely
\begin{align} \label{eq:closure_eq}
    \begin{aligned}
    &\rho_\ell(p) = \hat{\rho}_\ell e^{c_\ell (p - \hat{p})}\\
    &k_{r, \ell}(S_\ell) = S_\ell^2
    \end{aligned}
     & \quad \text{in } \Omega^T,
\end{align}
where $\hat{\rho}_\ell$ is a reference value for the density, in \sib{\kilogram\per\cubic\meter}, $c_\ell$ a phase specific compressibility, in \sib{\per\pascal}, and $\hat{p}$ a reference pressure value, in \sib{\pascal}. 
\end{subequations}

By using \eqref{eq:saturation_eq}, the problem \eqref{eq:mass_bal_equi} can be recast in an equivalent form in terms only 
of one saturation, here $S_0$, and pressure. Furthermore, we assume the porosity to be time independent and we replace one mass balance with the sum of the two mass balances. Thus, we replace \eqref{eq:mass_bal_equi_eq} and \eqref{eq:saturation_eq} with
\begin{subequations}\label{eq:equidimeq}
\begin{align}
    \begin{aligned}
    & \phi \partialder{[\rho_0 S_0 + \rho_1 (1-S_0)]}{t} + \divergence Q_T = f_T \\
    & \phi\partialder{(\rho_0 S_0)}{t} + \divergence Q_0 = f_0
    \end{aligned}
    \quad \text{in } \Omega^T,
\end{align}
where $f_T = f_0 + f_1$, and $Q_T$ is the total flux, defined as
\begin{gather}
    Q_T =  - \sum_{\ell=0,1}   \rho_\ell \lambda_\ell K \left(\nabla p + \rho_\ell g\nabla z \right)
    \quad \text{in } \Omega^T,
\end{gather}
\end{subequations}
%

\noindent
where $\lambda_1$ is now written as a function of $S_0$ instead of $S_1$. 
This formulation is useful to highlight the different nature of the variables, $p$ and $S_0$, with consequent advantages in the discretization scheme, see \Cref{sec:discretization}. Indeed, problem \eqref{eq:equidimeq} shows the so-called mixed parabolic-hyperbolic behaviour if the compressibility is taken into account, an elliptic-hyperbolic behaviour otherwise. As a consequence, the pressure varies smoothly while the saturation may be discontinuous in $\Omega$ \cite{Trangenstein1989, Trangenstein1989a}. 

\subsection{Dimensional analysis} \label{subsec:dimensional}

To describe the flow regime and compare the results of different simulations, we can scale the equations, initial and boundary conditions, and consequently obtain dimensionless groups. We denote by $x_\mathrm{ref}$ a reference value for the generic quantity $x$, the corresponding dimensionless variable is $\Tilde{x} = x / x_\mathrm{ref}$. At each physical variable or data we associate a reference value, some of them might depends on each other as a consequence of the Pi theorem \cite{Buckingham1914}. Neglecting the source term for simplicity, the mass balance becomes
\begin{gather}\label{eq:scaled_eq_1}
    \partialder{(\Tilde{\phi}\Tilde{\rho_\ell} S_\ell)}{t} + \frac{t_\mathrm{ref}K_\mathrm{ref}}{\phi_\mathrm{ref}L_\mathrm{ref}\mu_\mathrm{ref}} \Tilde{\nabla}\cdot \left[ \Tilde{\lambda_\ell}\Tilde{K}\left( \frac{p_\mathrm{ref}}{L_\mathrm{ref}}\Tilde{\nabla}\Tilde{p} + \rho_\mathrm{ref} g_\mathrm{ref} \Tilde{\rho_\ell}\Tilde{g}\Tilde{\nabla} \Tilde{z} \right)  \right] = 0.
\end{gather}
It governs the fluid motion in $\widetilde{\Omega^T} = \widetilde{\Omega}  \times \left(0, {T}/{t_\mathrm{ref}}\right]$, where $\widetilde{\Omega}$ is given by scaling each dimension of $\Omega$ by $L_{\mathrm{ref}}$, with analogous initial and boundary condition of \eqref{eq:bc_ic}.
We set $t_\mathrm{ref} = {\phi_\mathrm{ref}L_\mathrm{ref}} /{u_\mathrm{ref}}$ and retrieve a reference velocity from a finite Darcy-type equation $u_\mathrm{ref} = {K_\mathrm{ref} p_\mathrm{ref}} / ({\mu_\mathrm{ref}L_\mathrm{ref}})$, thus we set a reference pressure from viscous quantities $p_\mathrm{ref} = {\mu_\mathrm{ref}^2}/({K_\mathrm{ref}\rho_\mathrm{ref}})$. Replacing them into \eqref{eq:scaled_eq_1}, we obtain
\begin{equation}\label{eq:dimless_mass_bal}
    \partialder{(\Tilde{\phi}\Tilde{\rho_\ell} \Tilde{S_\ell})}{t} +  \Tilde{\nabla}\cdot \left[ \Tilde{\lambda_\ell}\Tilde{K}\left( \Tilde{\nabla}\Tilde{p} + E_A \Tilde{\rho_\ell}g\Tilde{\nabla} \Tilde{z} \right)  \right] = 0,
\end{equation}
where $E_A = {\phi_\mathrm{ref} \rho_\mathrm{ref}^2 g_\mathrm{ref} L_\mathrm{ref} K_\mathrm{ref}}/{\mu_\mathrm{ref}^2}$ is a dimensionless group that indicates the ratio of the effects of the gravity forces and the viscous forces. It can be seen as an Archimedes' number specific to two-phase flow problem driven by the gravity, or an adaptation of the gravity number \cite{Riaz2007, Fumagalli2012, Lee2015} to scenarios where the reference velocity is ambiguous, such as the case of countercurrent flow driven by gravity.


\subsection{A mixed-dimensional model}\label{sec:mixed-dimensional}
As mentioned before, when fractures are present we rely on dimensionally reduced models to approximate fractures with lower dimensional objects immersed in the rock domain. Consequently, we need to devise a new set of partial differential equation, derived from mass balance and Darcy law, that describe the flow in the fractures and the interaction with the rock matrix. 
The dimensional reduction by itself is well-established, and we refer to  \cite{Martin2005, Angot2009, Formaggia2014} for single phase flow and in \cite{Jaffre2011, Fumagalli2012d, Brenner2015, Ahmed2017, Brenner2018} for two-phase flow.
When this technique is combined with the representation of the geometry introduced in Section \ref{sec:domain}, we arrive at the following governing equations.



The flow problem is: find $(S_{0, i}, p_i)$ in each
$\Omega^T_i=\Omega_i\times (0,T]$, $i = 1,\ldots,I$ and $\zeta_{\ell,j}$ on each $\Gamma_j^T = \Gamma_j \times (0, T]$, $j=1,\ldots,J$, for both phases $\ell=0, 1$, such that
\begin{subequations}\label{eq:full_system}
    \begin{align}
            &\varepsilon_i^{a_i}\phi_i \partialder{[\rho_{0,i} S_{0,i} + \rho_{1,i}(1-S_{0,i})]}{t} + \varepsilon_i^{a_i} \divergence Q_{T,i} + \sum_\ell \sum_{\Gamma_j\in \hat{S}_i} \Xi_j^i \mathcal{N}_j(\rho_{\ell}\lambda_{\ell})\zeta_{\ell,j} = f_{T, i}, \label{eq:press_eq} \\
            &\varepsilon_i^{a_i}\phi_i \partialder{(\rho_{0,i} S_{0_i})}{t} + \varepsilon_i^{a_i} \divergence Q_{0,i} + \sum_{\Gamma_j\in\hat{S}_i}\Xi_j^i \mathcal{N}_j(\rho_{0} \lambda_{0})\zeta_{0,j} = f_{0,i},
            \label{eq:mass_bal_mix} 
    \end{align}
with the following constitutive law for the mortar fluxes, $\zeta_{\ell, j}$,
\begin{gather}\label{eq:claw_mortar}
    \zeta_{\ell,j} - \varepsilon_l^{b_j-1} k_{\perp, j} \left\{ \frac{2}{\varepsilon_l}[\Pi_j^h tr(p_h) - \Pi_j^l p_l] - \mathcal{N}_j(\rho_{\ell}) g\nabla z\cdot \upsilon_{out, j} \right\} = 0,
\end{gather}
and the following boundary conditions at all times to close the previous system:
\begin{align}
    \begin{aligned}
    &Q_\ell |_{\partial_j \Omega_i} \cdot \upsilon_{out, j} - \Xi_j^h\mathcal{N}_j(\rho_\ell \lambda_\ell)\zeta_{\ell,j} = 0, \quad & \forall j \in \check{S}_i  \\
    &Q_\ell \cdot \upsilon |_{\partial \Omega_i}= 0, \quad & \text{on} \ \partial \Omega_i \setminus \partial_j\Omega_i
    \end{aligned}
\end{align}
\end{subequations}
Finally, appropriate initial conditions for the primary variables $S_{0,i}$, $p_i$ and $\zeta_j$ have to be provided for all the domains $\Omega_i$ and mortars $\Gamma_j$

The divergence and gradient in \eqref{eq:press_eq} and \eqref{eq:mass_bal_mix} are meant as being tangential to the manifold associated to the considered domain. Since we only consider  flat subdomains, these operators can be easily written in local intrinsic orthogonal, and fixed, coordinates.

The additional terms in \eqref{eq:press_eq} and \eqref{eq:mass_bal_mix} compared to \eqref{eq:equidimeq} are linked to the mortar fluxes and describe the interactions between domains of different dimensions.


where we recall that with the indeces $l$ and $h$ we denote $\varepsilon_l = \varepsilon_i, \Omega_i \in \check{R}_j$, $p_h = p_i$, $\Omega_i \in \hat{R}_j$, and $p_l = p_i$, $\Omega_i \in \check{R}_j $, $tr$ is a trace operator, $k_{\perp, j}$ the normal permeability. Equation \eqref{eq:claw_mortar} derives from the flux described in \eqref{eq:mass_bal_equi_eq} along the orthogonal direction of the mortar $\Gamma_j$, where the pressure gradient is approximated by a finite difference across the subdomains. 
In \eqref{eq:claw_mortar},
$\zeta_\ell$ can be seen as a volumetric flux or velocity 
divided by the mobility, it is used also to compute the interface upwind direction as described below.
 

We call the formulation of the problem \eqref{eq:full_system} dual-mortar formulation, where we name \eqref{eq:press_eq} pressure equation and \eqref{eq:mass_bal_mix} mass balance.

\subsection{Manipulation of the model for discretization}
The mass fluxes, $Q_{\ell,i}$ and $Q_{T,i} = Q_{0,i} + Q_{1,i}$, can be expressed differently, leading to different discretizaton methods, as we show in \Cref{sec:discretization}.
The first straightforward option is to write $Q_{\ell,i} = \rho_{\ell,i} q_{\ell,i}$. Consequently,
\begin{equation}\label{eq:total_flux}
    Q_{T,i} = \rho_{0,i} q_{0,i} + \rho_{1,i} q_{1,i},    
\end{equation}
where $q_{\ell,i}$ is the volumetric flux given by the Darcy law:
\begin{equation}\label{eq:volumetric_fluxes}
    q_{\ell,i} = -K_i\lambda_{\ell,i}(\nabla p_i + \rho_{\ell,i} g\nabla z), \quad \ell = 0,1    
\end{equation}
This formulation is adopted in the phase-potential upwind (PPU) discretization, detailed in \Cref{sec:discretization_ppu}.
Note that the physical properties of a fracture subdomain are now described by two quantities, the normal permeability, $k_{\perp,j}$ and the in-plane permeability $K_i$, that govern, respectively, the flow across and along the subdomain.

The second option separates the physical contribution of the fluxes, highlighting the one due to the pressure gradient, called viscous flux, and the one due to gravity. In this case, the total flux in the pressure equation \eqref{eq:press_eq} is still written as the sum of the phases contribution, as in \eqref{eq:total_flux}, while the flux $Q_0$ in the mass balance equation \eqref{eq:mass_bal_mix} is further subdivided as

\begin{equation}\label{eq:flux_mass_bal}
    Q_{0,i} =  V_{0,i} + G_{0,i},
\end{equation}
being $V_{0,i}$ the viscous mass flux and $G_{0,i}$ is the flux in mass due to gravity, as a result of the different densities of the fluids.
\begin{align}
    &V_{0,i} = \rho_{0,i} \frac{\lambda_{0,i}}{\lambda_{T,i}} (q_{0,i} + q_{1,i})\label{eq:V}\\
    &G_{0,i} = \rho_{0,i} K_i\frac{\lambda_{0,i} \lambda_{1,i}}{\lambda_{T,i}}(\rho_{1,i} - \rho_{0,i})g\nabla z\label{eq:G}
\end{align}
with $\lambda_{T,i} = \lambda_{0,i} + \lambda_{1,i}$ the total mobility, in \sib{\pascal\per\second}. This second formulation is adopted in the hybrid upwind (HU) discretization, detailed in \Cref{sec:discretization_hu}. 


At the discrete level, we need that the operator $\mathcal{N}_j$ that maps variables from domains $\Omega_i$ to an adjacent the mortar domain $\Gamma_j$ has an upwinding nature, to avoid the appearance of spurious oscillations. To this purpose, we use the following formulation:
%
\begin{align}\label{eq:interface_upwind}
    \mathcal{N}_j(u_{\ell}) = 
    \begin{cases}
        \Pi_j^i u_{\ell,i} \ & \text{if} \ \zeta_{\ell, j} < 0, \\ 
        \Pi_j^i tr(u_{\ell,i}) \ & \text{if} \ \zeta_{\ell, j} \geq 0, 
    \end{cases}
\end{align}
%
where $u_{\ell,i}$ is defined in $\Omega_i \in\check{R}_j$ in the first case, and $u_{\ell,i}$ is defined in $\Omega_i\in\hat{R}_j$ in the second case.
Because of the absence of internal fluxes, zero-dimensional domains should be treated carefully.
For each zero-dimensional $\Omega_i$, for $k=1,\ldots,I$ and for both phases $\ell=0,1$, the mass balances reduce to finding $(S_{0, i}, p_i)$ such that
\begin{align}\label{eq:full_system_0d}
    \begin{aligned}
        &\varepsilon^{a_i}\phi_i \partialder{[\rho_{0,i} S_{0,i} + \rho_{1,i}(1-S_{0,i})]}{t} +  \sum_\ell \sum_{\Gamma_j\in \hat{S}_i} \Xi_j^i \mathcal{N}_j(\rho_{\ell} \lambda_{\ell})\zeta_{\ell,j} = f_{T, i},\\ 
        &\varepsilon^{a_i}\phi_i \partialder{(\rho_{0,i} S_{0,i})}{t} + \sum_{\Gamma_j\in\hat{S}_i}\Xi_j^i \mathcal{N}_j(\rho_{0} \lambda_{0})\zeta_{0,j} = f_{0,i}.
    \end{aligned}
\end{align}
See Fig.~\ref{fig:mixed_dim_domain}(b) for an example of three branches intersecting.

We model the fluid mobility with a quadratic dependence on the saturation:
\begin{gather}\label{eq:claw_lambda}
    \lambda_\ell = S_\ell^2.    
\end{gather}
%



\section{Discretization}\label{sec:discretization}

In this section, we present the considered approach to solve system \eqref{eq:full_system} numerically. The main numerical challenge is how to approximate the
non-linear advective part appropriately since it is composed by multiple terms that might act quite differently.
This has an impact also on the convergence of the Newton method chosen to solve the resulting non-linear problem.


After a brief introduction to the discretization tools in \Cref{subsec:common}, we will present the first discretization method
in \Cref{sec:discretization_ppu}, which will be referred as PPU. The second discretization
method, which is one of the original contribution of this work, is an extension to the
mixed-dimensional framework of the hybrid upwind scheme and its application to the case
of fractured media \cite{Lee2015, Hamon2017, Hamon2016, Hamon2016a, Hamon2018, Hamon2020, Bosma2022}.
It is presented in \Cref{sec:discretization_hu} and it is referred as HU. The aim of this method is address the countercurrent flow problem more effectively,  reducing
the number of Newton iteration required to solve the non-linear system at each timestep.

\subsection{Common setup}\label{subsec:common}

\begin{figure}[h!]
    \centering
    \includegraphics[width=0.3\textwidth]{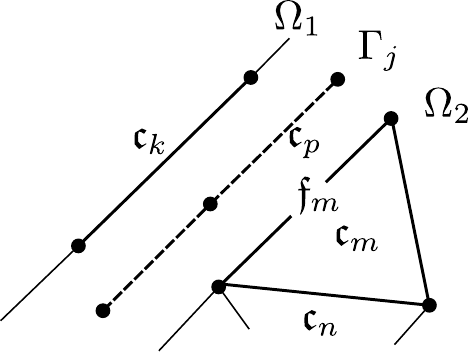}
    \caption{Grid elements of a 1D domain, $\Omega_1$ facing a 2D domain, $\Omega_2$, and mortar,
        $\Gamma$. Any cell is denoted by $\mathfrak{c}$, while the faces, i.e., the boundary of
        $\mathfrak{c}$ are $\mathfrak{f}$.}
    \label{fig:grids_el}
\end{figure}

The subdomains $\Omega_i$ and mortars $\Gamma_j$ are approximated with grids composed by
simplex or hexahedron cells $\mathfrak{c}$ of characteristic size $h_\mathfrak{c}$ and Lebesgue measure $|\mathfrak{c}|$. The
faces of a cell $\mathfrak{c}$ are denoted by $\mathfrak{f}$. Each face between two cells
$\mathfrak{c}_m$ and $\mathfrak{c}_n$ is associated with a unit normal $\hat{\upsilon}_{mn}$. We
denote by $adj(\mathfrak{c}_m)$ the neighbouring cells of $\mathfrak{c}_m$. The
formulation with mortar fluxes let us easily deal with non-conforming meshes across the interfaces between domains of different dimensions, as depicted in Fig.~\ref{fig:grids_el}.

We solve \eqref{eq:full_system} by discretizing the
equations with a cell centred finite volume method \cite{Blazek2015, Barth2017}. 


In the following, we indicate with $x_{mn}$ the evaluation of
a generic variable, $x$, on the face $\mathfrak{f}$ between cells $\mathfrak{c}_m$ and $\mathfrak{c}_n$. Moreover, we denote the jump of a variable $x$ between two cells $\mathfrak{c}_m$ and $\mathfrak{c}_n$ as $\Delta x_{mn} = x_m - x_n$. Before continuing, we define the discrete upwind operator $\mathcal{M}$ between $\mathfrak{c}_m$ and $\mathfrak{c}_n$ with respect to an upstream direction
$v_{mn}$ and the oriented segment $\overline{mn}$ that connects the cell centers from
$\mathfrak{c}_m$ to $\mathfrak{c}_n$:
\begin{gather*}
    \mathcal{M}(x_m, x_n; v_{mn}) =
    \begin{dcases}
        x_m & \qif\  v_{mn} \cdot \overline{mn} \geq 0, \\
        x_n & \qif\ v_{mn} \cdot \overline{mn} < 0.
    \end{dcases}
\end{gather*}

Similarly, the component of a variable $\m x$ projected on a cell $\mathfrak{c}_p$ of the mortar $\Gamma_j$ with the interface upwind defined in \eqref{eq:interface_upwind}, according to the upstream direction $v_p$ defined on the mortar cell $\mathfrak{c}_p$, is computed as:
%
%
\begin{gather*}
    [\mathcal{N}_j(\m x; v_p)]_p = \begin{dcases}
        \Pi_{j,pk}^i x_k & \qif \ v_p < 0, \ x_k \in \Omega_i \in \check{R}_j,                         \\
        \Pi_{j,pn}^i x_n & \qif \ v_p \geq 0, \ x_n \in \Omega_i \in \hat{R}_j, \ n \in S_j^{\mathfrak{c}},
    \end{dcases}
\end{gather*}
where $S_j^{\mathfrak{c}}$ is the set of indices $n$ denoting the cells of $\Omega_i$ at the
boundaries facing $\Gamma_j$. To ease  notation, form now on we denote $[\mathcal{N}_j(\m x; v_p)]_p$ as $x_p^{v}$. 


A remark should be made for the discrete maps, $\Pi$ and $\Xi$, in particular
when a non-matching discretization is used between objects of different dimensions. 
To achieve maximum accuracy, in the case of cell-wise constant variables as are applied herein, the map for intensive variables (e.g. the pressure) should apply area-weighted averaging, while extensive quantities (e.g. fluxes) should be summed. 




Let us consider a cell $\mathfrak{c}_p \in \Gamma_j$, a cell   $\mathfrak{c}_k \in \Omega_i \in \check{R}_j$ and a boundary face $\mathfrak{f}_m$ of cell $\mathfrak{c}_m \in \hat{R}_j$, see Fig.~\ref{fig:grids_el}. 
Let $l_{mp}$ and $l_{kp}$ be the overlap surface between boundary face $\mathfrak{f}_m$ and cell $\mathfrak{c}_p$, and cell $\mathfrak{c}_k$ and cell $\mathfrak{c}_p$, respectively.
For the intensive variable we have: $\Pi_{j,mp}^{h,ave} = {l_{mp}}/{|\mathfrak{c}_p|}$ and $\Pi_{j,kp}^{l,ave} = {l_{kp}}/{|\mathfrak{c}_p|}$,
while for an extensive variable: $\Pi_{j,mp}^{h,sum} = {l_{mp}}/{|\mathfrak{f}_m|}$ and $\Pi_{j,kp}^{l,sum} = {l_{kp}}/{|\mathfrak{c}_k|}$. Similarly
for $\Xi$. To ease the notation, we will not specify the type of map, it should be clear from the context. 
As mentioned before, it is worth noting that the faces and cells facing mortars are, in general, non-conforming, as depicted in Fig.\ref{fig:grids_el}.

%

\subsection{Standard upwinding}\label{sec:discretization_ppu}
Let $S_i^{\mathfrak{f}}$ be the set of indexes for faces $\mathfrak{f}$ on the boundary towards
$\Gamma_j$ and let $\circ$ be the element-wise (Hadamard) product. 
The discrete primary variables, i.e. the degree of freedom (d.o.f.) associated to the cell centers, are
the components of the vectors $\m p_i$, $\m S_{0,i}$ for the subdomains $\Omega_i$ and
$\m\zeta_{0,j}$, $\m\zeta_{1,j}$ at the mortars $\Gamma_j$.
The discretization in space
of \eqref{eq:full_system} with fluxes from \eqref{eq:total_flux} and \eqref{eq:volumetric_fluxes}, reads for $\ell=0,1$
\begin{gather}\label{eq:discrete_full_system_ppu}
    \begin{aligned}
        &\varepsilon_i^{a_i} \partialder{(\m u_{0,i} + \m u_{1,i})}{t} + \varepsilon_i^{a_i} \mm{D} \, \m Q_{T,i}^{PPU} + \m \psi_i^p  = 0,                                                                                 && \text{in } \Omega_i, \\
        &\varepsilon_i^{a_i}\partialder{\m u_{0,i}}{t} + \varepsilon_i^{a_i} \mm{D} \, \m Q_0^{PPU} + \m \psi_i^s = 0,                                                                                                     && \text{in } \Omega_i, \\
        &\m\zeta_{\ell,j} - \varepsilon_l^{b_j-1}k_{\perp,j}|\underline{\mathfrak{c}}_j|\circ \left[ \frac{1}{\varepsilon/2}(\mm\Pi_j^h \underline{\text{tr}}(\m p_h) - \mm \Pi_j^l \m p_l) - \m \rho_\ell^{\zeta_{\ell,j}} g\nabla z\cdot \upsilon_{out, j} \right] = 0, && \text{in } \Gamma_j\ .
    \end{aligned}
\end{gather}
%
The discrete boundary conditions, for each vector element $k$, are given by:
\begin{gather*}
    \begin{aligned}
         & Q_{\ell,i,k}^{PPU} - \Xi_{j,km}^i\rho_{\ell,m}^{\zeta_{\ell,j}}\lambda_{\ell,m}^{\zeta_{\ell,j}}\zeta_{\ell,j,m} = 0, & k \in S_i^{\mathfrak{f}}, \\
         & Q_{\ell,i,k}^{PPU} = 0,                                                            & \text{otherwise.}
    \end{aligned}
\end{gather*}
%
 
To ease the notation, in what follows, we drop the subscript denoting subdomain and interfaces and let relations hold for each subdomain $\Omega_i$ or mortar $\Gamma_j$.

The terms of the accumulation variables are $u_{0,m} = \phi_m\rho_{0,m} S_{0,m}$, and
$u_{1,m} = \phi_m \rho_{1,m}(1- S_{0,m})$. The elements of discrete divergence operator
$\mm D$ are defined as:
\begin{gather*}
    D_{mn} = \frac{1}{|\mathfrak{c}_m|}
    \begin{dcases}
        1  & \qif \ \mathfrak{c}_n \in adj(\mathfrak{c}_m) \ \text{and} \ \hat{\upsilon}_{mn} \ \text{points inwards} \ \mathfrak{c}_m,  \\
        -1 & \qif \ \mathfrak{c}_n \in adj(\mathfrak{c}_m) \ \text{and} \ \hat{\upsilon}_{mn} \ \text{points outwards} \ \mathfrak{c}_m, \\
        0 & \qif \ \mathfrak{c}_n \notin adj(\mathfrak{c}_m).
    \end{dcases}
\end{gather*}
The integrated-normal mass flux
$Q_{\ell,mn}^{PPU} = \rho_{\ell,mn}^{PPU} q_{\ell,mn}^{PPU}$
is discretized adopting an upwind scheme whose upstream direction is linked to the phase
potential, $\Phi_\ell$, \eqref{eq:phase_potential}. In particular, we have: $\rho_{\ell,mn}^{PPU} = \mathcal{M}(\rho_{\ell,m},\rho_{\ell,n}; \Tilde{q}_{\ell,mn})$, where $\Tilde{q}_{\ell,mn}$ represents the phase potential multiplied by the transmissibilities. It is computed as
\begin{equation*}
    \underline{\Tilde{q}}_{\ell} = \mm T \, \m p + \mm T^g g \m\rho_\ell \circ \m z, 
\end{equation*}
where $\mm T$ and $\mm T^g$ are the transmissibilities, that depend on the permeability $K_a$ and the grid
element geometry. Note that the elements of $\underline{\Tilde{q}}_\ell$ are $\Tilde{q}_{\ell,mn}$, i.e., quantities evaluated at the faces. In the present work, $\mm T$ is computed with the multi-point flux
approximation (MPFA) \cite{Aavatsmark2002, Agelas2008}, although other methods could be used, such as the two point flux approximation
(TPFA) \cite{Aavatsmark2002, Aavatsmark2007a}. 
$\mm T^g$ is the consistent trasmissibility for the gravity term, its determination is explained in \cite{Starnoni2019}. 
The integrated volumetric flux is:
$q_{\ell,mn}^{PPU} = \lambda_{\ell,mn}^{PPU} \Tilde{q}_{\ell,mn}$,
where $\lambda_{\ell, mn}^{PPU}$  are phase mobility defined as
$\lambda_{\ell, mn}^{PPU} = \mathcal{M}(\lambda_{\ell,m}, \lambda_{\ell,n}; \Tilde{q}_{\ell,mn})$.
The elements of the source term, $\m \psi^p$ in the pressure equation and $\m \psi^s$ in the
mass balance for the phase $0$ are respectively given by
\begin{gather*}
    \begin{aligned}
         & \psi_m^p = \sum_\ell \sum_{\Gamma_j\in \check{S}} \Xi_{j,mn} (\rho_{\ell,n}^{\zeta_{\ell,j}} \lambda_{\ell,n}^{\zeta_{\ell,j}} \zeta_{\ell,j,n}) - f_{\ell,m} , \\
         & \psi_{m}^s = \sum_{\Gamma_j\in\check{S}} \Xi_{j,mn} (\rho_{0,n}^{\zeta_{0,j}} \lambda_{0,n}^{\zeta_{0,j}} \zeta_{0,j,n}) - f_{0,m}.
    \end{aligned}
\end{gather*}
%
Regarding the constitutive law associated with the mortar fluxes, we highlight that it is not a
partial differential equation but an algebraic expression, so no spatial discretization scheme is required.

\subsection{Hybrid upwind}\label{sec:discretization_hu}

The semi-discrete in space counterpart of the problem \eqref{eq:full_system}, with fluxes from \eqref{eq:total_flux}, \eqref{eq:volumetric_fluxes} for the pressure equation \eqref{eq:press_eq} and \eqref{eq:flux_mass_bal} for the mass balance equation \eqref{eq:mass_bal_mix}, is given by for $\ell = 0,1$
\begin{gather}\label{eq:discrete_full_system_hu}
    \begin{aligned}
        &\varepsilon_i^{a_i}\partialder{(\m u_{0,i} + \m u_{1,i})}{t} + \varepsilon_i^{a_i} \mm{D} \ \m Q_{T,i} + \m \psi_i^p  = 0,                                                                             && \text{in } \Omega_i,  \\
        &\varepsilon_i^{a_i}\partialder{\m u_{0,i}}{t} + \varepsilon_i^{a_i} \mm{D} \left( \m V_{0,i} + \m G_{0,i} \right) + \m \psi_i^s = 0,                                                                   && \text{in }  \Omega_i, \\
        &\m\zeta_{\ell,j} - \varepsilon_l^{b_j-1} k_{\perp,j} \left[ \frac{1}{\varepsilon/2}(\mm\Pi_j^h \underline{\text{tr}}(\m p_h) - \mm \Pi_j^l \m p_l) -  \m\rho_\ell^{\zeta_{\ell,j}} g\nabla z\cdot \upsilon_{out, j} \right] = 0, && \text{in } \Gamma_j.
    \end{aligned}
\end{gather}
Concerning the different nature of the pressure equation \eqref{eq:press_eq} and mass balance \eqref{eq:mass_bal_mix}, they are discretized
adopting different strategies. The variables defining $\m Q_T$ are discretized with a blended
method that smoothly switches from a centred scheme to an upwind scheme, according to the
intensity of the jump of the phase potential evaluated at each face. Conversely, the fluxes in the mass balance are 
discretized with a pure upwind
scheme as described below. In particular, by dropping the domain index, we have
\begin{gather}\label{eq:Q_T}
    Q_{T,mn} = \sum_{\ell = 0,1}  \rho_{\ell, mn}  q_{\ell,mn}^{WA},
\end{gather}
where the densities are averaged at the faces according to the saturation:
\begin{gather*}
    \rho_{\ell,mn} = \frac{S_{\ell,m}\rho_{\ell,m}+S_{\ell,n}\rho_{\ell,n}}{S_{\ell,m}+S_{\ell,n}},
\end{gather*}
and the discrete volumetric fluxes integrated along the faces, $q_{\ell,mn}^{WA}$, are computed
following \cite{Bosma2022}. For the reader convenience, we report here the main steps. The
discrete volumetric normal flux, $q_{\ell,mn}^{WA}$, incorporates the grid and rock properties, the
fluid properties, and it is proportional to the phase potential:
\begin{gather}\label{eq:q_ell}
    q_{\ell,mn}^{WA} = \lambda_{mn}^{WA} T_{mn} \Delta \Phi_{\ell,mn},
\end{gather}
where, in the present work, the transmissibilities $T_{mn}$ are computed with the TPFA, although other methods could be used. The jump of the phase potential is defined as
\begin{equation*}
    \Delta \Phi_{\ell,mn} = \Delta p_{\ell,mn} + \rho_{\ell,mn}g\Delta z_{mn}.
\end{equation*}
The mobilities $\lambda_{mn}^{WA}$ are the core of the blended centred-upwind scheme: they are
the weighted average (WA) dependent
solution of the mobilities computed at the cells:
\begin{gather*}
    \lambda_{\ell,mn}^{WA} = \beta_{\ell,mn} \lambda_{\ell,m} + (1-\beta_{\ell,mn})\lambda_{\ell,n}
\end{gather*}
where the weight $\beta_{\ell,mn} \in [0,1]$ in an increasing function of the jump of the phase potential
$\Delta \Phi_{\ell, mn}$:
\begin{equation*}
    \beta_{\ell,mn} = 0.5 + \frac{1}{\pi}\arctan(c_{\ell,mn} \Delta \Phi_{\ell, mn}),
\end{equation*}
with $c_{\ell,mn}$ a coefficient that depends on the grid and fluid properties
\begin{equation*}
    c_{\ell,mn} = \min\left( \frac{ (k_{r,\ell}(1))^{-1} \max_{S_\ell}|k_{r,\ell}''(S_\ell)|}{\rho_{\ell,mn}},10^6 \right).
\end{equation*}
We remind that $k_{r, \ell}$ is the relative permeability of phase $\ell$. Details about $\beta_{\ell,mn}$ can be found in \cite{Bosma2022}.

The fluxes in the mass balance are discretized with an upwind scheme to honor the hyperbolicity
of the equation. Since the motion of the fluid is forced by two main physical driving forces,
one related to the pressure gradient and the other linked to the gravity field, these two
quantities are treated differently to upwind the variables.
The discretization of the viscous flux is a one-sided scheme that considers the total velocity
as upwind direction:
\begin{gather*}
    V_{0,mn} = \rho_{0,mn}^V\frac{\lambda_{0,mn}^V}{\lambda_{T,mn}^V} q_{T,mn}
\end{gather*}
where $\rho_{0,mn}^V = \mathcal{M}(\rho_{0,m}, \rho_{0,n}; q_{T,mn})$ and similarly for
$\lambda_{0,mn}^V$ and $\lambda_{T,mn}^V$. The total volumetric flux (or total velocity),
$q_{T,mn} = q_{0,mn} + q_{1,mn}$, is computed from \eqref{eq:q_ell}.

The gravity flux $G_{0,mn}$ is computed with a one-sided scheme with the upwind direction
dependent on the gravity effects. In particular, we have
\begin{gather*}
    G_{0,mn} = \rho_{0,mn}^G q_{0,mn}^G,
\end{gather*}
where the upstream direction for the density is the volumetric flux due to the gravity effects,
$\rho_{0,mn}^G = \mathcal{M}(\rho_{0,m}, \rho_{0,m}; q_{0,mn}^G)$. The volumetric flux,
$q_{0,mn}^G$, is computed as:
\begin{gather}\label{eq:qG}
    q_{0,mn}^G = T_{mn} \frac{\lambda_{0,mn}^G \lambda_{1,mn}^G}{\lambda_{T,mn}^G}(\rho_{0,mn} - \rho_{1,mn})g\Delta z_{mn}
\end{gather}
where $\lambda_{0,mn}^G = \mathcal{M}(\lambda_{0,m}, \lambda_{0,n}; \omega_{0,mn})$ and
analogously for the other mobilities. The function $\omega_{\ell,mn}$ describes the gravity
effects acting on fluids with different densities:
\begin{gather*}
    \omega_{\ell,mn} = \lambda_{k,mn}^g \left( (\rho_{k,mn}-\rho_{\ell,mn})g\Delta z_{mn} \right),
\end{gather*}
where the upwinded mobility is given by
\begin{gather*}
    \lambda_{k,mn}^g = \begin{dcases}
        \lambda_{k,m} & \qif \ \rho_{k,mn} < \rho_{\ell,mn}, \\
        \lambda_{k,n} & \qif \ \rho_{k,mn} \geq \rho_{\ell,mn}.
    \end{dcases}
\end{gather*}
Further details can be found in \cite{Bosma2022}.

The constitutive laws of the mortar fluxes, \eqref{eq:claw_mortar}, do not evidence the different
behaviour of the pressure and saturation variables. As in the PPU approach, no particular
discretization method are required.
Future developments regard the study of a pressure-mass formulation of the mortar fluxes
constitutive laws and a suitable discretization method.

The spatial order of convergence of the discretization of the fluxes $Q_T$, $V_0$, and $G_0$ is at most one for smooth solutions \cite{Hamon2016}. The discretization of the mortar fluxes is expected to be of order one since an error proportional to $h_\mathfrak{c}$ is introduced
with the approximation of $\mathcal{N}_j$. This occurs because no trace reconstruction is used to retrieve the value of $x$ on
$\partial_j\Omega_i$, instead, the value at the cell center is used. Thus, we expect the spatial order of convergence of the discretization scheme to be at most one for smooth solutions. A numerical verification is performed in \Cref{sec:case_1}. 

\subsection{Time discretization and numerical details}\label{sec:discretization_misc}

For the time marching, we use the implicit Euler scheme for its well known stability properties
that allow us to use large timesteps. An implicit method leads to a system of discrete non-linear
equations to be solved at each time step, which is solved with the Newton method. Note that the
Newton method requires the knowledge of the Jacobian, which changes, at each iteration.

Several methods can be adopted for the computation of the Jacobian or its approximation, such
as a finte difference, complex step \cite{Squire1998, Lantoine2012} or automatic
differentiation \cite{Bischof2008}. We adopt the latter, which is exact up to machine precision.

In case of convergence failure of Newton method, the time step is halved and the iterative method is restarted to the previous step. The procedure is
repeated till convergence is reached or the time step becomes excessively small for a practical application and it is thus
stopped due to lack of convergence. We set the threshold for the timestep to be $10^{-12}$.

The discretization method needs to comply with the constraint on the saturation, that is $S_\ell \in [0,1]$. We force its enforcement by clipping the saturation value after each Newton iteration.

We implement the discretization methods in PorePy, a simulation tool for fractured and
deformable porous media suited for the mixed-dimensional problem \cite{Keilegavlen2020}.


\section{Numerical validation}\label{sec:numerical_validation}

We consider three test cases to evaluate the performance of the numerical methods described in
the previous sections. The first one, presented in \Cref{sec:case_1}, is a 2D domain with one
single fracture. It is divided into three sub-cases in which the fracture orientation and rock
properties are modified. The second case, reported in \Cref{sec:case_2}, involves a fracture
network composed of 10 fractures with different permeability, some of which intersecting.
In the last test case, shown in
\Cref{sec:case_3}, we study a 3D geometry cut by 8 fractures with multiple intersections.

To keep the presentation simple, we apply the same fluid properties to the two phases for all
the test cases: the heavy phase, denoted by the subscript $0$, has a density $\rho_0 = 1$ and
dynamic viscosity $\mu_0 = 1$, the light phase is described by $\rho_1 = 0.5$ and $\mu_1 = 1$.

In all the cases, unless differently specified, we use the following convention to define the
dimensionless groups $E_A$ as in \eqref{eq:dimless_mass_bal}: the reference values are related
to the matrix domain where $L_\mathrm{ref}$ is the vertical length and $\rho_\mathrm{ref}$ is the difference
of the densities of the two phases.

The primary criteria to compare the discussed numerical schemes is the number of Newton
iterations done, being a proxy of the associated computational cost, and they accuracy to compute
the numerical solution. Other specific results are described case by case.

The iterative method is stopped when the error, $err$, defined as the norm of the normalized
$\ell_2$-norm of the increment, $\delta x$, is smaller than a given tolerance:
\begin{gather*}
    err = \frac{\|\delta x\|_2}{\sqrt{\dim{\delta x}}} < tol.
\end{gather*}
Unless differently specified, we set tolerance equal to $tol = 10^{-6}$.


\subsection{Case 1. Single fracture}\label{sec:case_1}

\begin{figure}[h]
    \centering
    \includegraphics[width=0.2\textwidth]{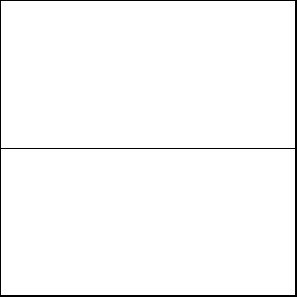}
    \hspace{0.02\textwidth}
    \includegraphics[width=0.2\textwidth]{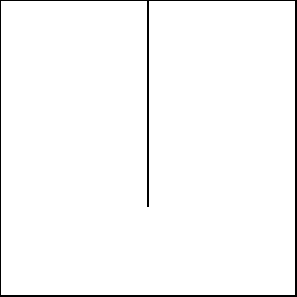}
    \hspace{0.02\textwidth}
    \includegraphics[width=0.2\textwidth]{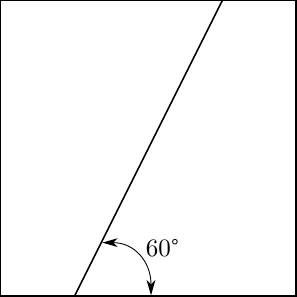}\\
    \subfloat[\textit{Case 1.a}]{\includegraphics[width=0.2\textwidth]{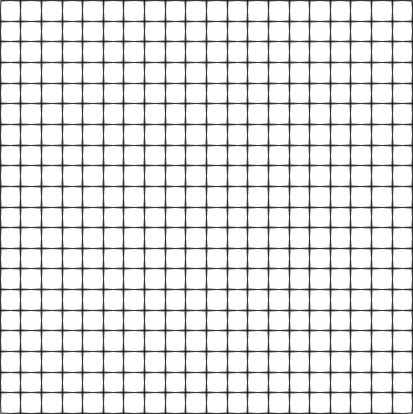}}
    \hspace{0.02\textwidth}
    \hspace{0.05cm}
    \subfloat[\textit{Case 1.b}]{\includegraphics[width=0.2\textwidth]{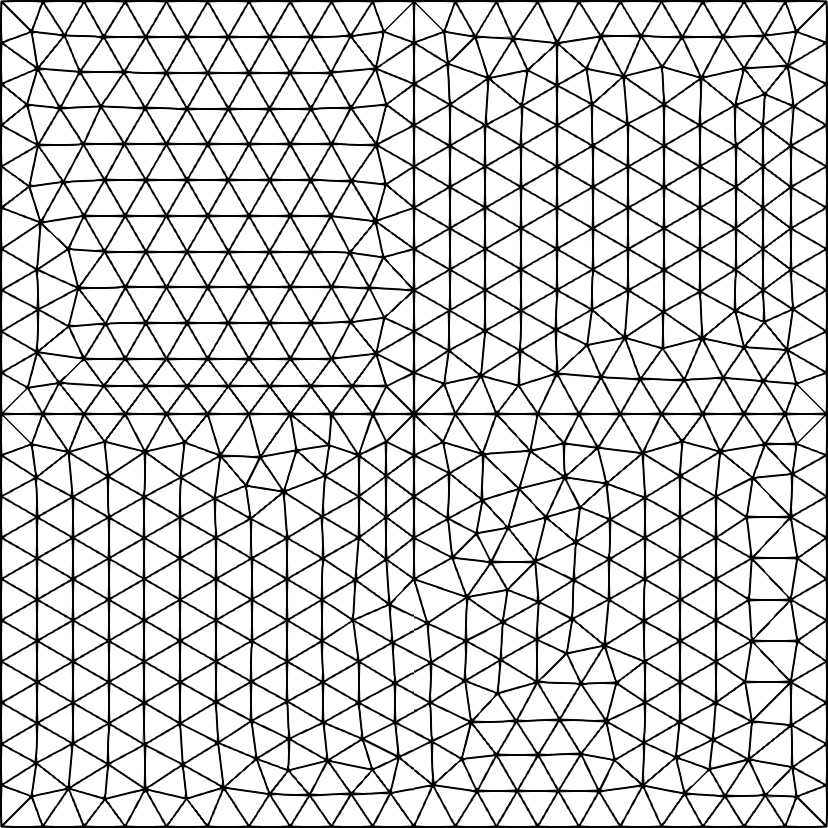}}
    \hspace{0.02\textwidth}
    \hspace{0.05cm}
    \subfloat[\textit{Case 1.c}]{\includegraphics[width=0.2\textwidth]{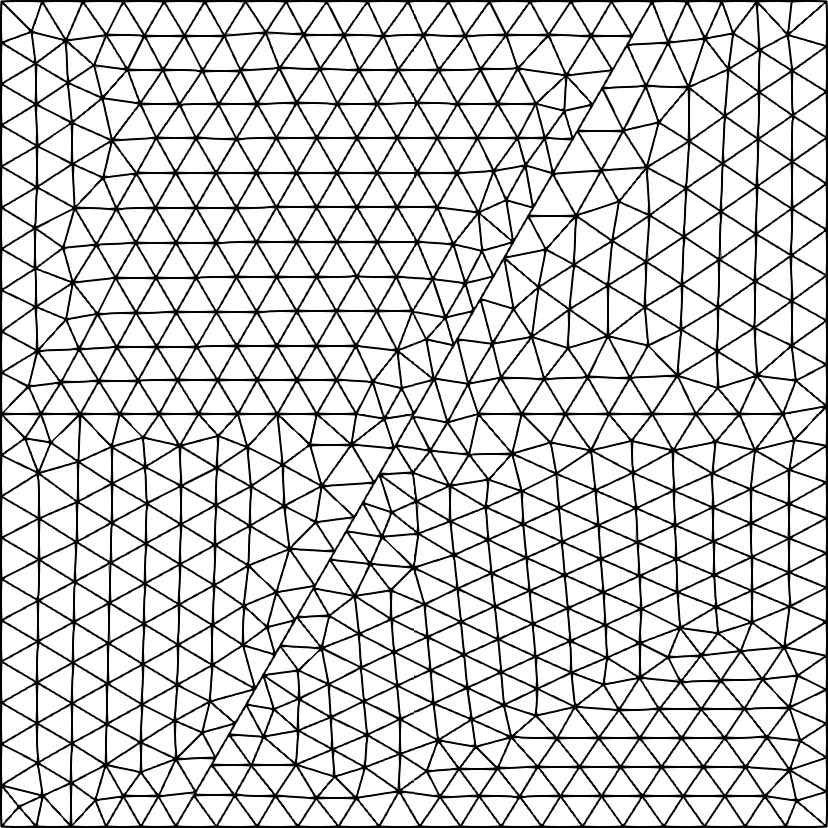}}
    \caption{Case 1. Domains and grids for the example in \Cref{sec:case_1}. The fracture is
        located in three different positions, (a) horizontally with tips touching the borders, (b)
        vertically with one tip touching the border and the other immersed in the matrix, (c) with
        an angle and the tips touching the borders. Different type of grids are used: (a) structured
        with quadrilateral elements, (b) unstructured with triangular elements, (c) unstructured
        non-conforming at the fracture interface with triangular elements.}
    \label{fig:case_1_domain}
\end{figure}

The first test considers a simple 2D unit square domain, Fig.~\ref{fig:case_1_domain}, and aims
to study the main characteristics of the discretization methods. We investigate three sub-cases
by varying the position of the fracture, its permeability, and the grid elements shape
(triangles or squares), see Fig.~\ref{fig:case_1_domain}. In particular, \textit{Case 1.a} has an horizontal
fracture touching both borders, 
\textit{Case 1.b} has a vertical highly permeable fracture that ends inside the matrix, the last one
\textit{Case 1.c} has an oblique fracture.

The dynamic is defined by the unstable initial condition with the heavy phase on the top and
the light one on the bottom, separated by a sharp front and forced by the gravity, as depicted,
for example, in  Fig.~\ref{fig:case_1_slanted_dyn}.

\clearpage

\subsubsection{Case 1.a. Horizontal fracture}

The parameters and properties used in this case are summarized in Tab.~\ref{tab:case_1_hor_prop}.
The flow motion is strongly dependent on the permeability of the fracture that
crosses the entire domain, since the two halves of the domain only interact through the mortar
fluxes.

The computational grid is made of quadrilateral elements, as depicted in Fig.~\ref{fig:case_1_domain},
with a limited number of elements, whose side size is $0.05$.

\begin{figure}[h]
    \centering
    %
    {
    \subfloat[Cumulative number of flips of the upwind directions.]{
        \begin{minipage}[b]{0.33\textwidth}
            \centering
            \includegraphics[width=1\textwidth]{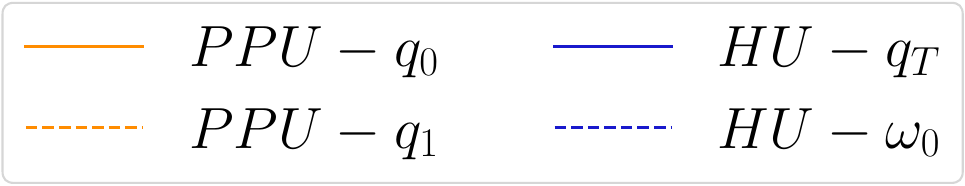}\\
            \begin{turn}{90}
                \hspace{1.5cm} \textbf{(1.a)}
            \end{turn}\includegraphics[scale=0.25]{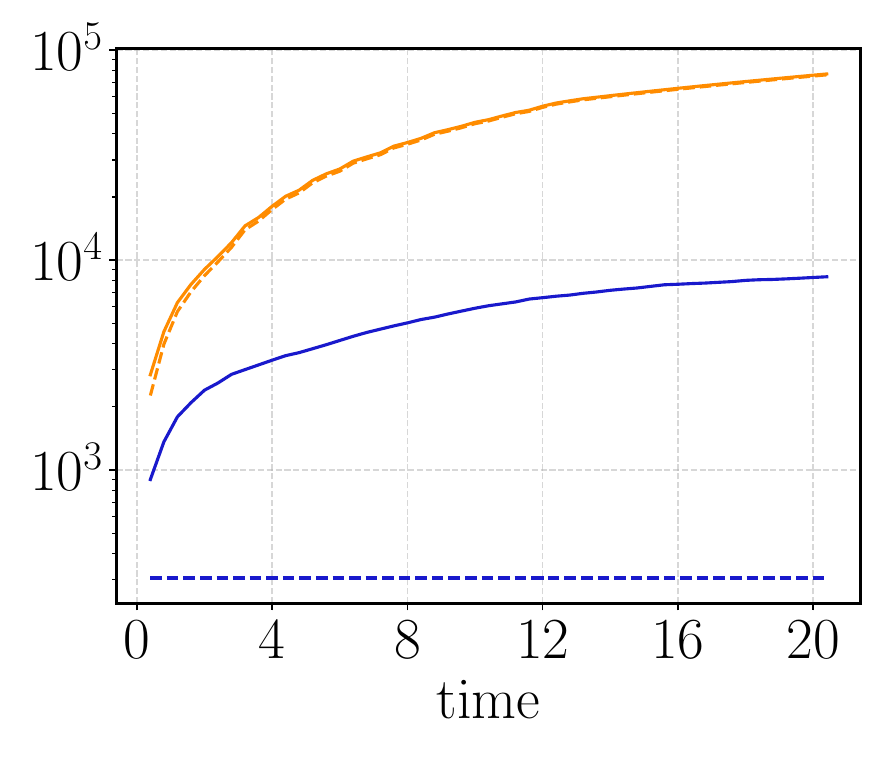}\\
            \begin{turn}{90}
                \hspace{1.5cm} \textbf{(1.b)}
            \end{turn}\includegraphics[scale=0.25]{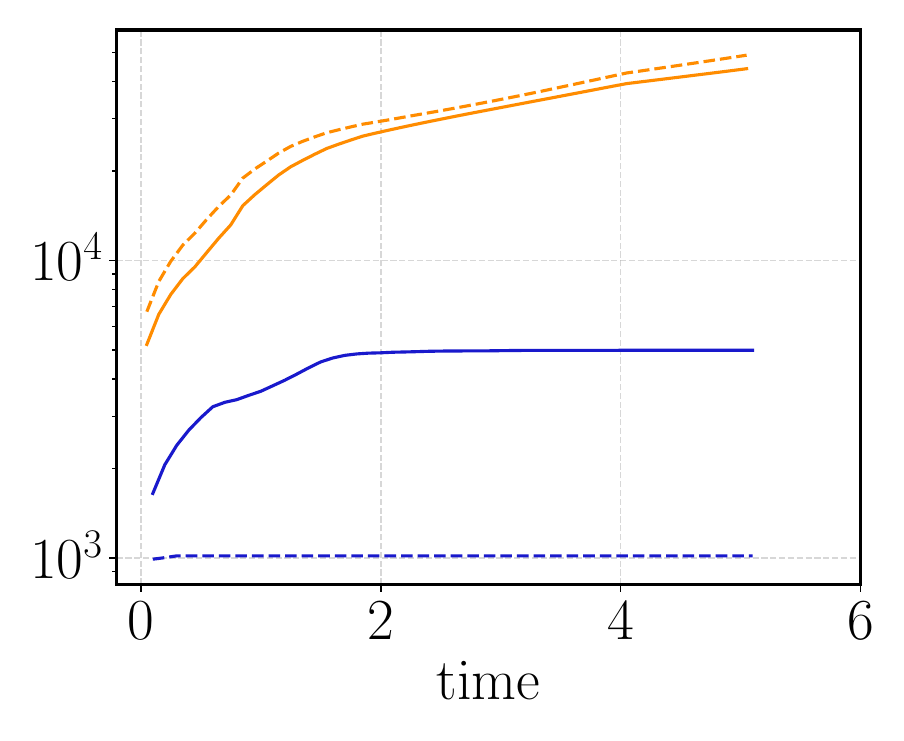}\\
            \begin{turn}{90}
                \hspace{1.5cm} \textbf{(1.c)}
            \end{turn}\includegraphics[scale=0.25]{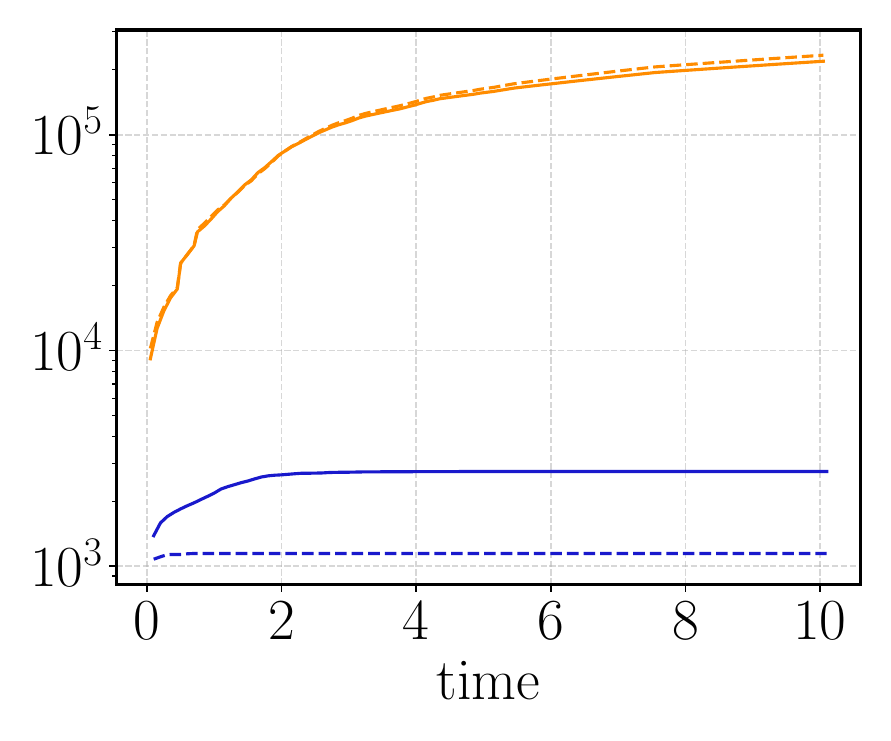}
            \\
            \begin{turn}{90}
                \hspace{1cm} \textbf{(1.c hc)}
            \end{turn}\includegraphics[scale=0.25]{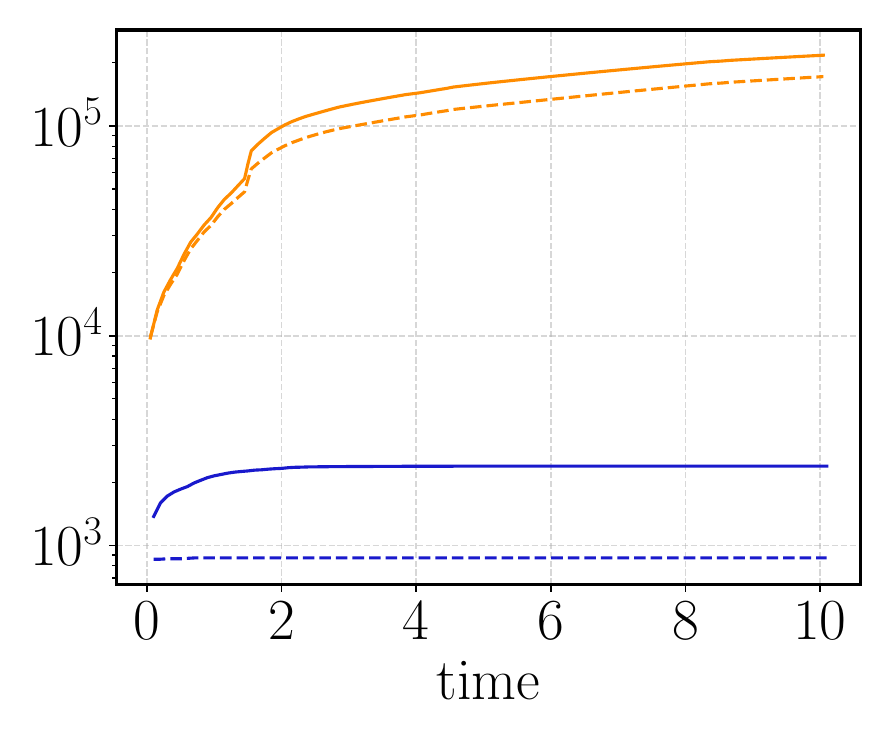}
        \end{minipage}
    }
    \subfloat[Cumulative number of time cuts.]{
        \begin{minipage}[b]{0.33\textwidth}
            \centering
            \includegraphics[width=0.75\textwidth]{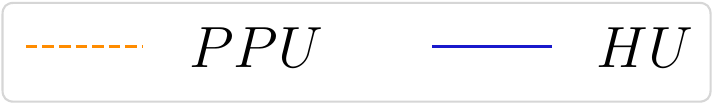}\\
            \includegraphics[scale=0.25]{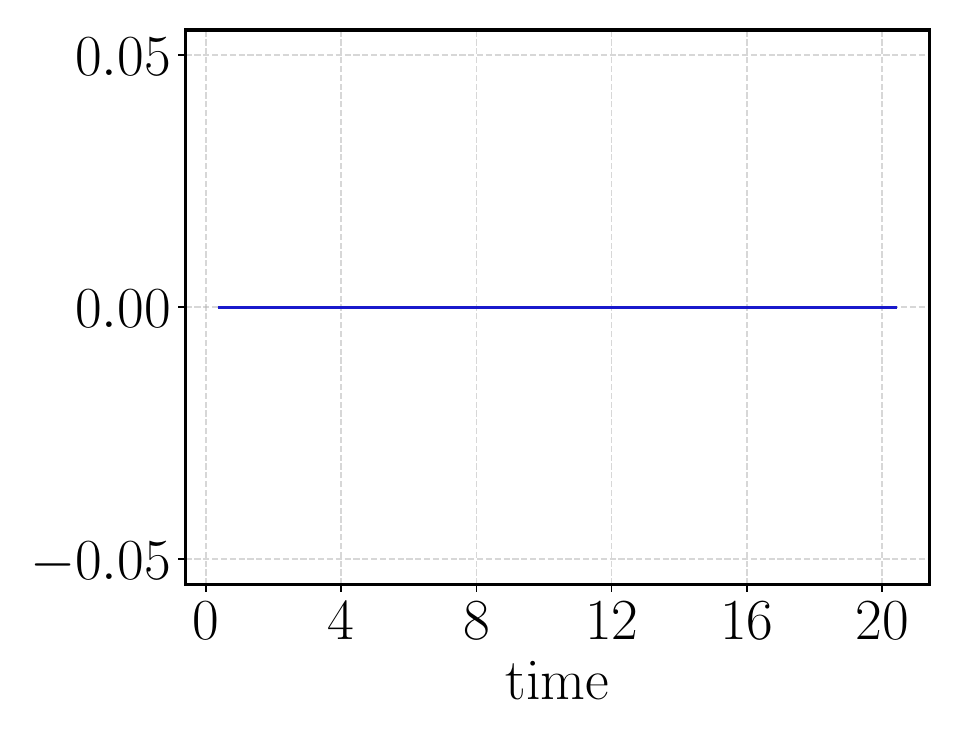}\\
            \includegraphics[scale=0.25]{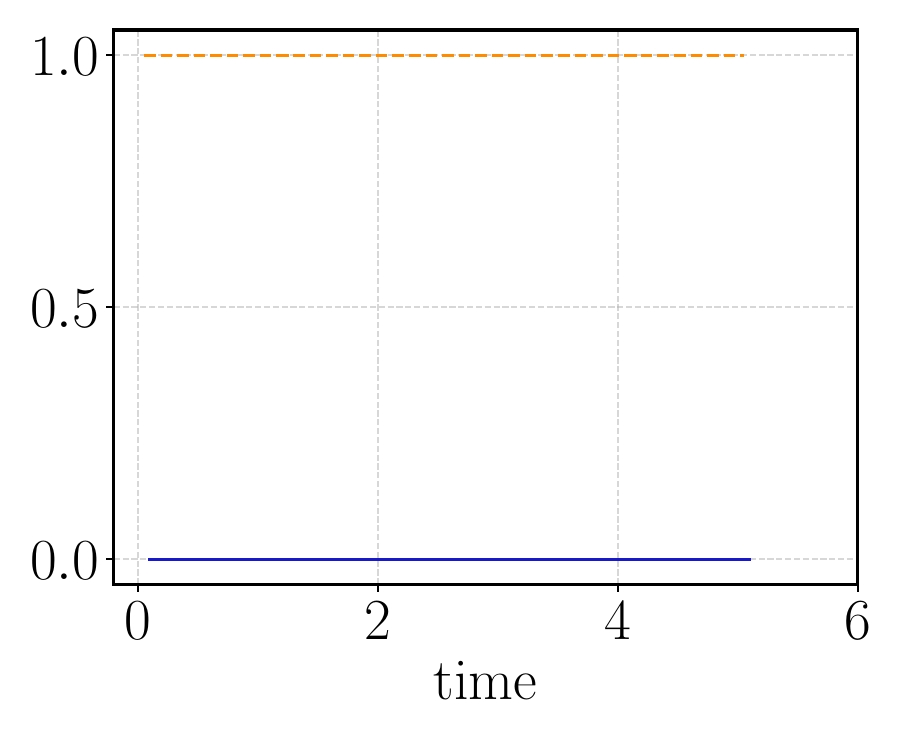}\\
            \includegraphics[scale=0.25]{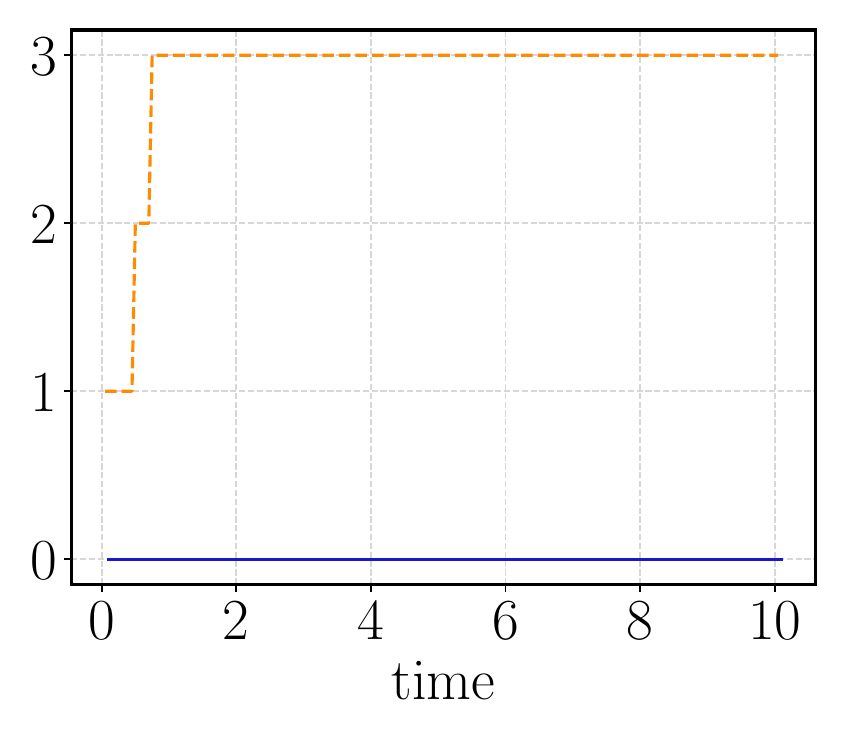}
            \includegraphics[scale=0.25]{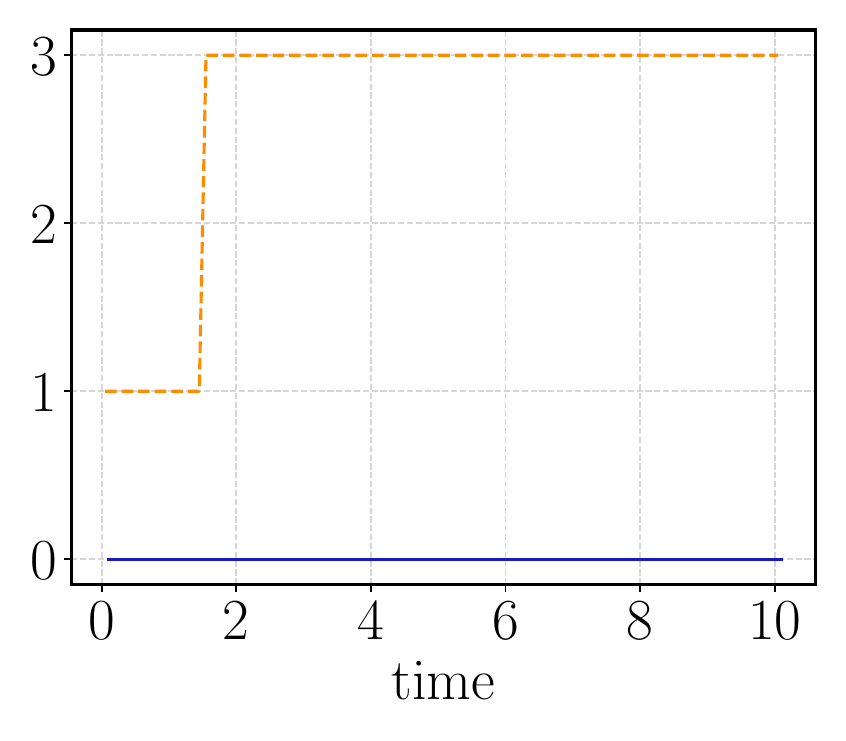}
        \end{minipage}
    }
    \subfloat[Cumulative number of Newton iterations.]
    {
        \begin{minipage}[b]{0.33\textwidth}
            \centering
            \includegraphics[width=1.1\textwidth]{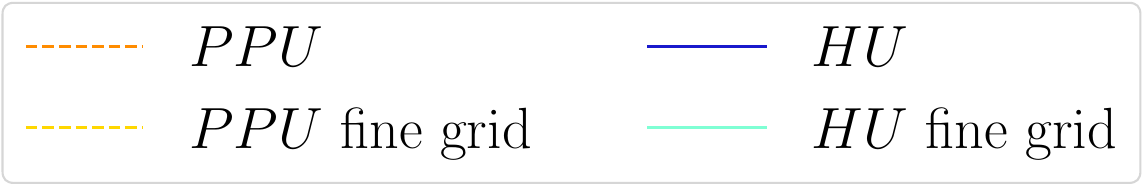}\\
            \includegraphics[scale=0.25]{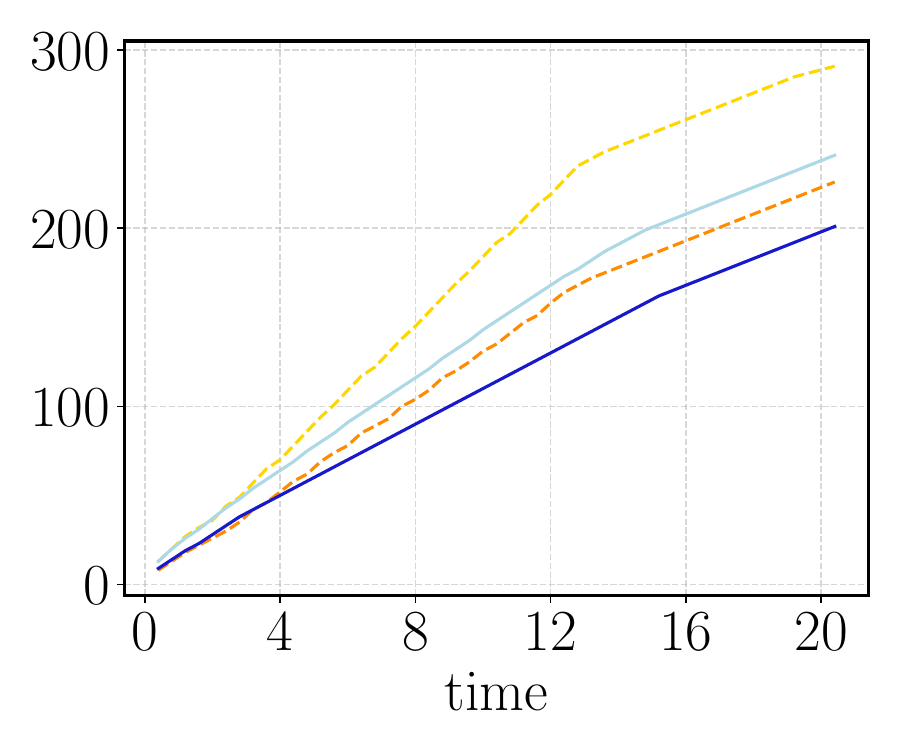}\\
            \includegraphics[scale=0.25]{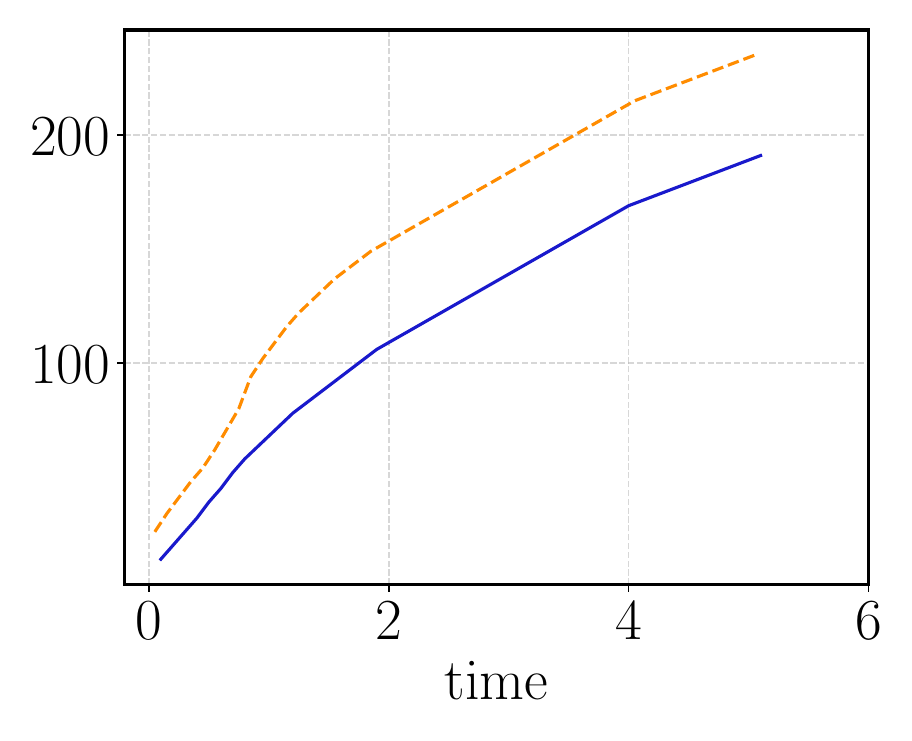}\\
            \includegraphics[scale=0.25]{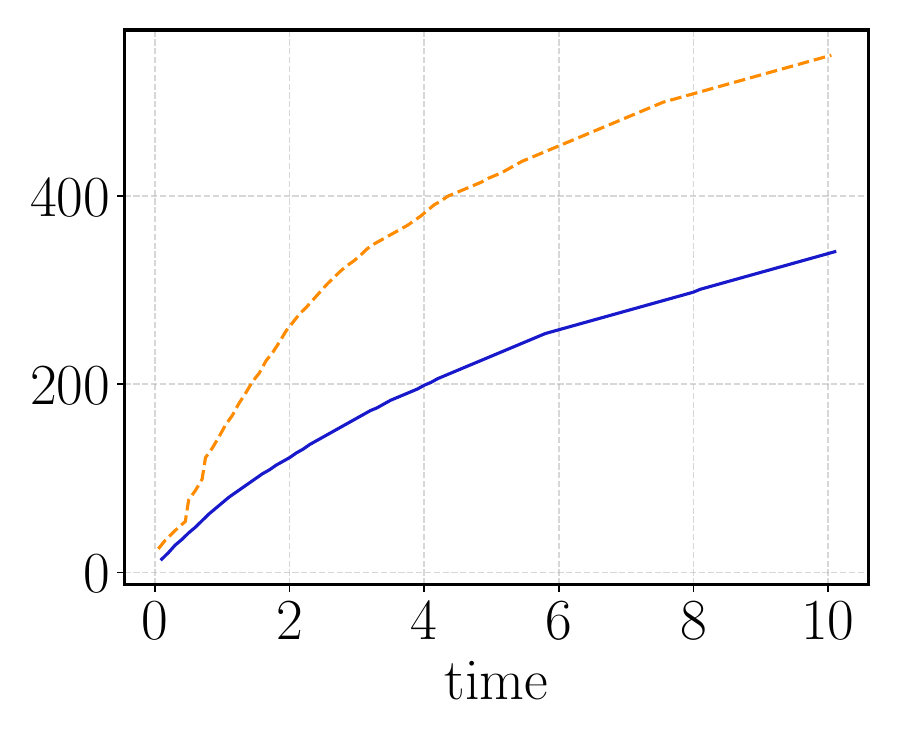}
            \includegraphics[scale=0.25]{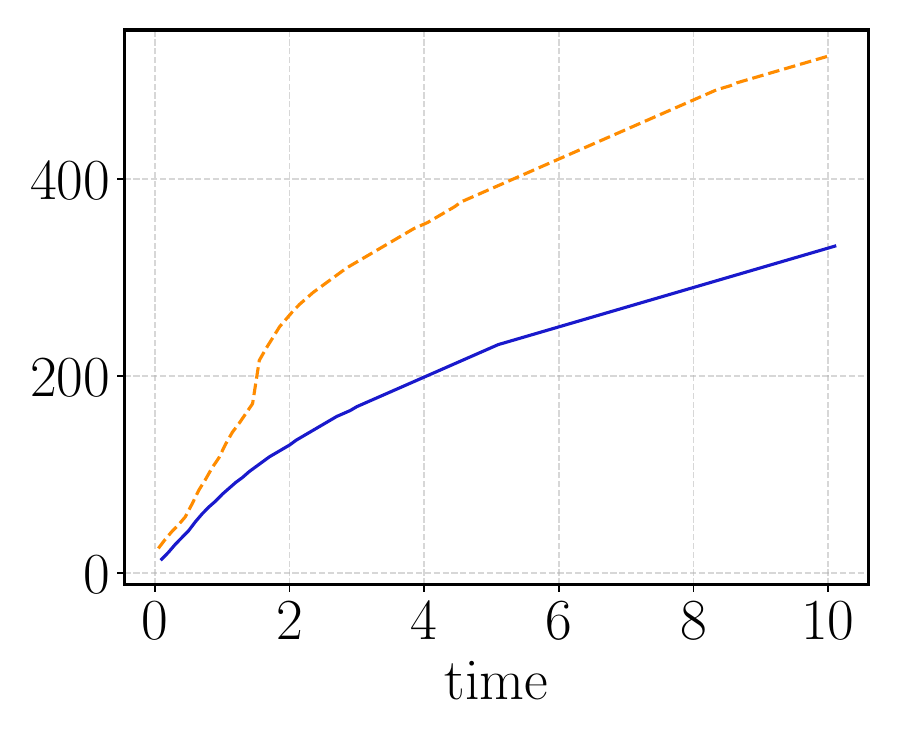}
        \end{minipage}
    } 
    }
    \\
    %
    {
    \subfloat[Cumulative number of flips of the upwind directions for \textit{Case 1.c}. 2D domain on the left panel, 1D domain at the center, mortar on the right.]{
    \begin{minipage}[b]{1.\textwidth}
        \includegraphics[scale=0.21]{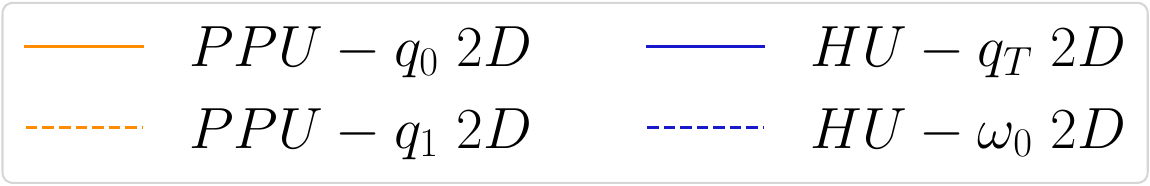}
        \includegraphics[scale=0.21]{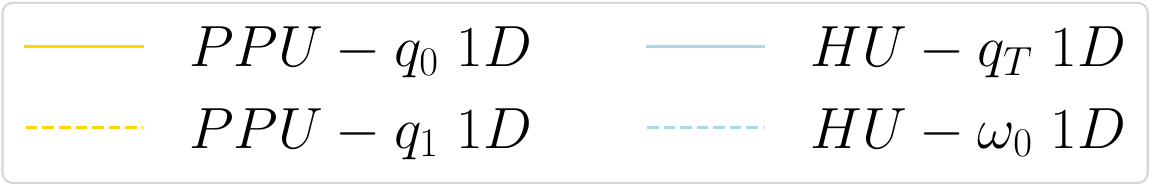}
        \includegraphics[scale=0.21]{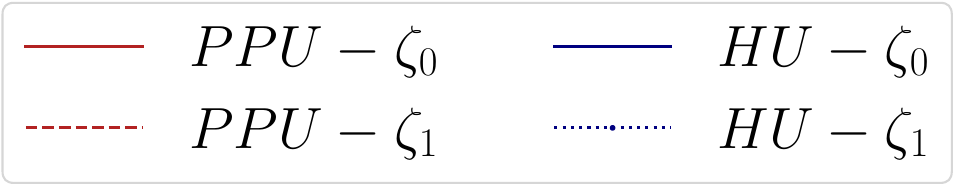}
        \begin{turn}{90}
            \hspace{1.5cm} \textbf{(1.c)}
        \end{turn}\includegraphics[scale=0.25]{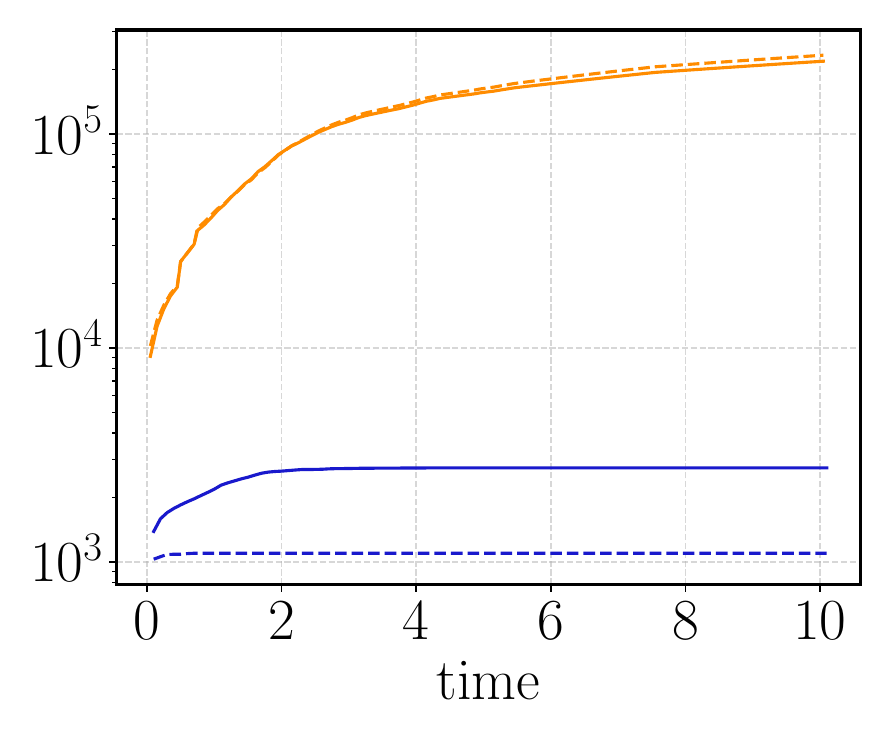}
        \includegraphics[scale=0.25]{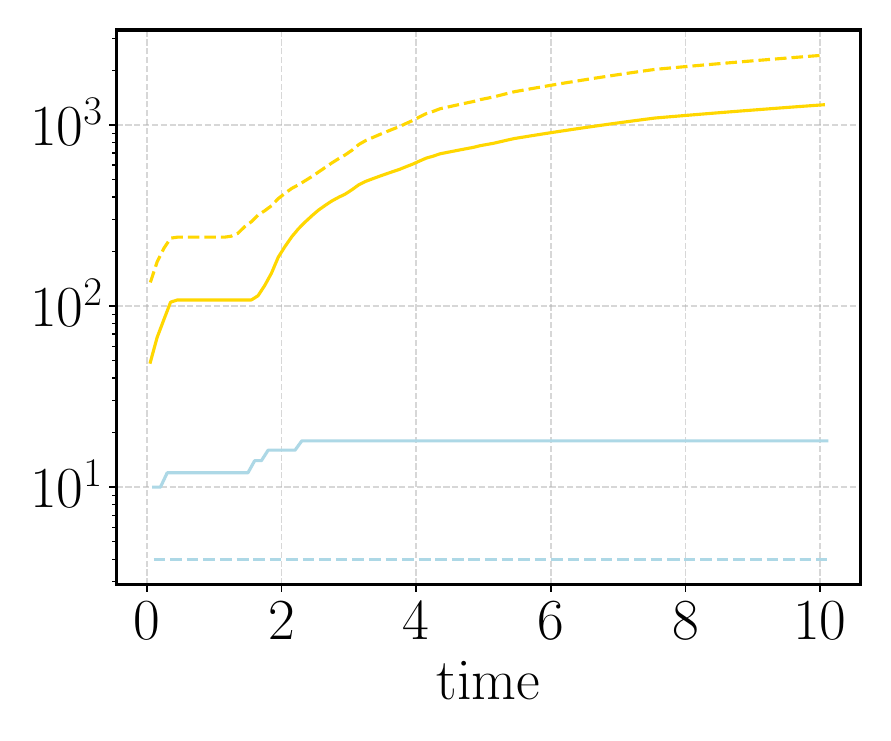}
        \includegraphics[scale=0.25]{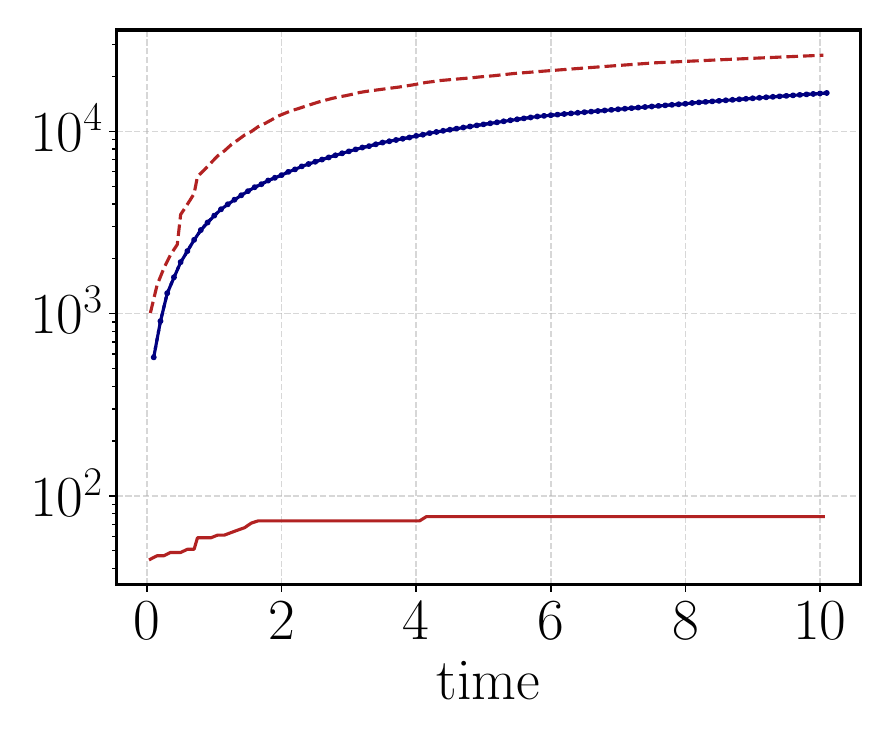}
    \end{minipage}
    }
    }
    \caption{Case 1. Characteristics regarding the iterative method.}
    \label{fig:case_1_newton}
\end{figure}

Due to the simplicity of the test case here a minor advantage is provided by the HU regarding
the cumulative number of Newton iterations, Fig.~\ref{fig:case_1_newton}. In the same picture,
we can see the time-cumulative number of wasted flips of the upwind direction, summed over all the grid
faces of the domain. Unless differently specified, only the 2D domain is shown, since it is the one that affects the Newton iteration the most. In particular, for the PPU the upwind directions are linked to the
volumetric fluxes, $q_0$ and $q_1$, that are more prone to change during the Newton iterations.
For the HU the directions are determined by the total velocity, $q_T$, and the gravity related
function $\omega_0$, that is very stable throughout the simulation. The middle panel of
Fig.~\ref{fig:case_1_newton} shows the cumulative number of timestep reductions, see
\Cref{sec:discretization_misc}. For this particular case, no cuts are performed. Regarding the
cumulative number of Newton iterations, we notice a small gap between the two methods, which
grows, and also the absolute values does, alongside the refinement of the grid. In
Fig.\ref{fig:case_1_newton} the cell size is reduced to $0.025$.

\begin{figure}
    \centering
    \includegraphics[width=0.3\textwidth]{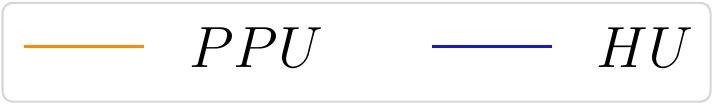}

    \subfloat[$t = 6.8$]{\includegraphics[width=0.4\textwidth]{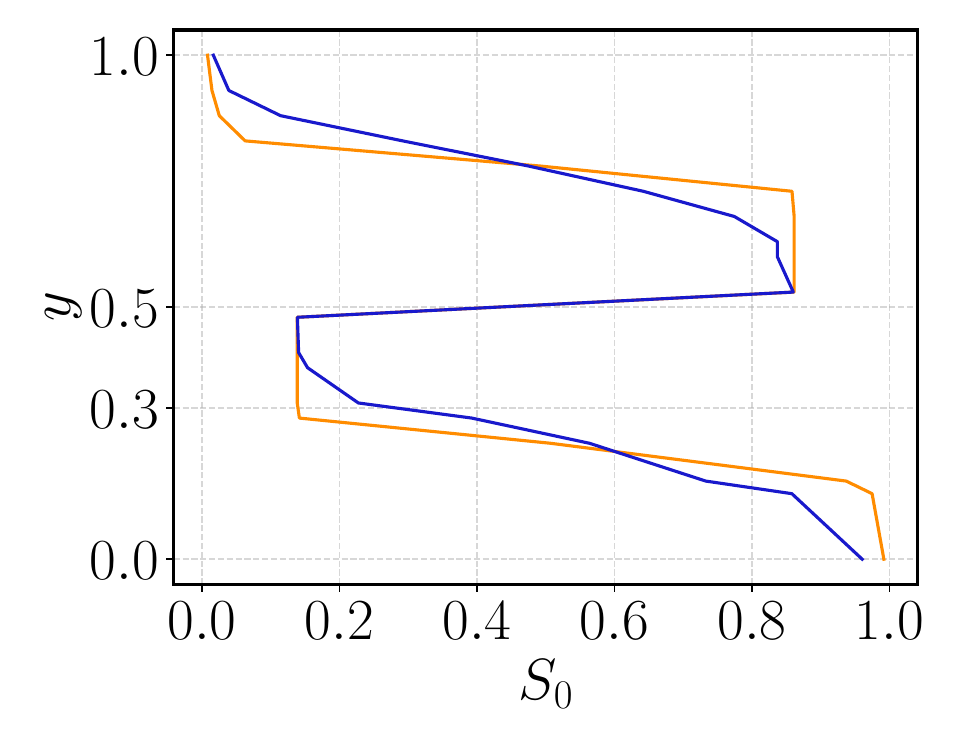}}
    \subfloat[$t = 20$]{\includegraphics[width=0.4\textwidth]{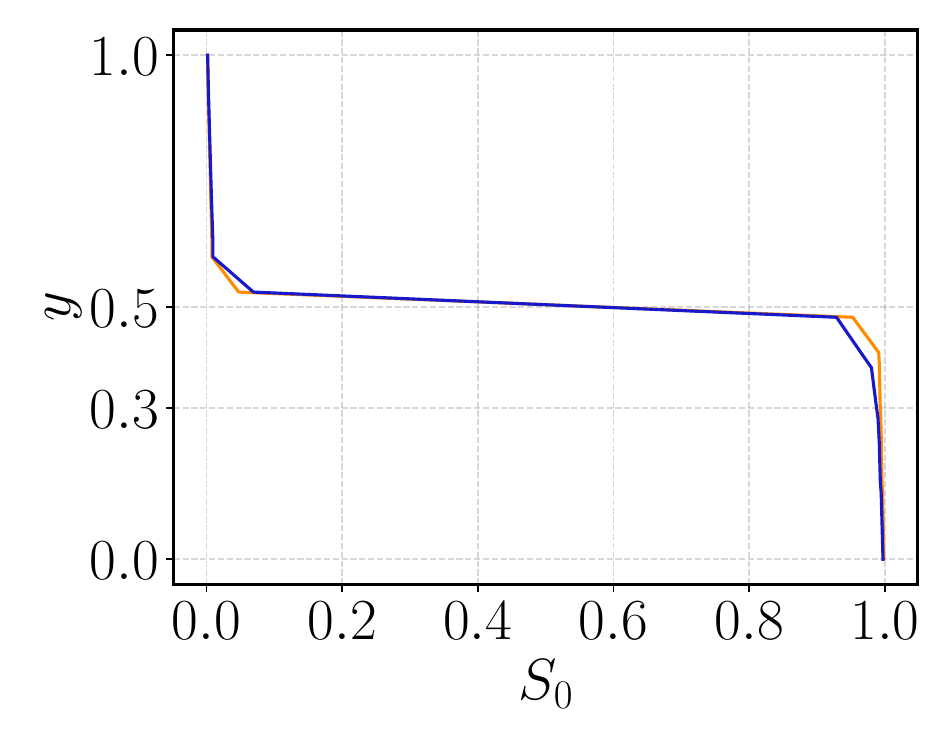}}

    \caption{Case 1. Horizontal. Saturation profile along a vertical line. The left panel illustrates
        $S_0$ at time $t = 6.8$, during the countercurrent flow of the two phases. A more diffusive trend
        is observable for the HU scheme than the PPU. In the right panel, showing the stationary solution
        at $t = 20$, the lines nearly coincide. No visible differences in the numerical diffusion are
        appreciable at the discontinuity since the interface fluxes are discretized with the same scheme.}
    \label{fig:case_1_hor_saturation}
\end{figure}

Note that the two methods behave differently in the presence of a sharp saturation front.
In particular, the HU method is more diffusive than the PPU, as shown in
Fig.~\ref{fig:case_1_hor_saturation}(a) that represents the saturation profile along a vertical
line during the transient. A sharp jump in the saturation is visible at $y=0.5$ because of the
presence of the fracture. The behavior nearby the fracture is similar for both PPU and HU since
the fluxes at the interface boundaries are computed with the same discretization scheme, see
\Cref{sec:discretization}.
\begin{table}[h!]
    \centering
    \begin{tabular}{ll}
        \hline
        Matrix intrinsic permeability   & $K_2 = 1$              \\
        Fracture intrinsic permeability & $K_1 = 1$              \\
        Fracture normal permeability    & $k_{\perp,1} = 0.1$    \\
        Fracture cross-sectional area   & $\varepsilon_1 = 0.01$ \\
        Matrix porosity                 & $\Phi_2 = 0.25$        \\
        Fracture porosity               & $\Phi_1 = 0.25$        \\
        Total simulation time           & $ t_{end} = 20$        \\
        Timestep max                    & $\Delta t_{max} = 0.4$ \\
        $E_A$                           & $0.0625$               \\
        \hline
    \end{tabular}
    \caption{Case 1. Horizontal. Parameters used for this problem.}
    \label{tab:case_1_hor_prop}
\end{table}
%



%

\clearpage

\subsubsection{Case 1.b. Vertical fracture}
The parameters and properties defining in this case are summarized in Tab.~\ref{tab:case_1_vertical_prop}.

\begin{figure}
    \centering
    \includegraphics[width=0.4\textwidth]{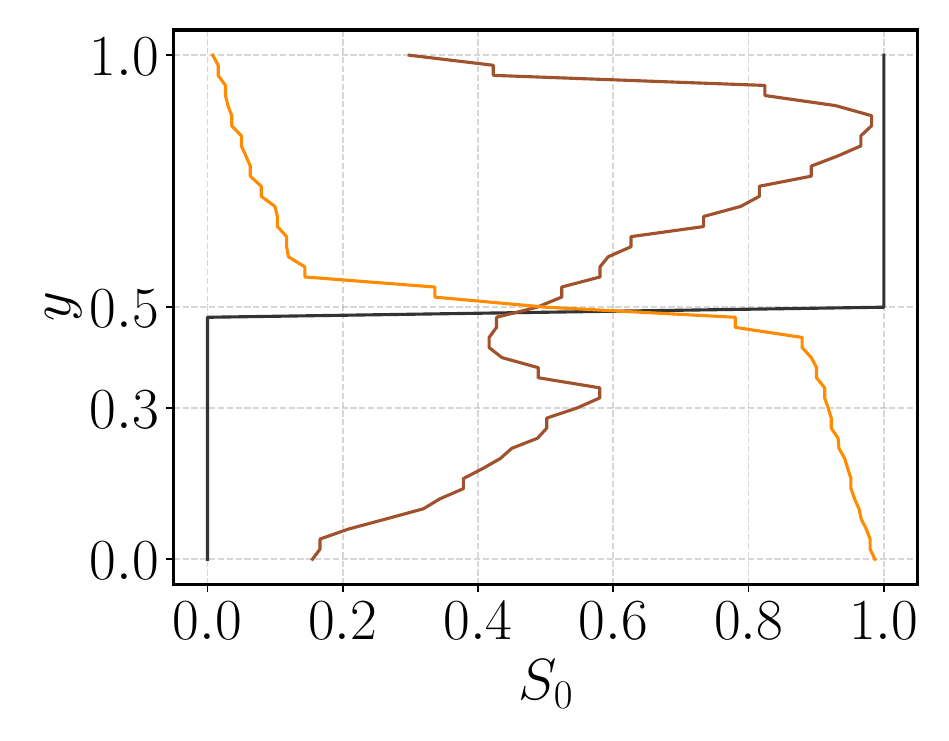} \\
    \includegraphics[width=0.4\textwidth]{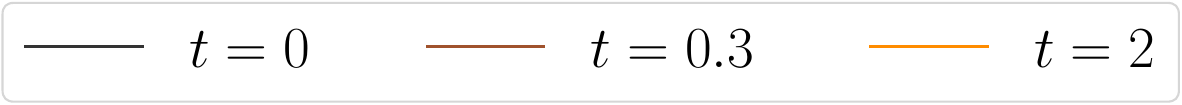}
    \caption{Case 1. Vertical. Saturation at different timesteps along a vertical line at the center of the 2D domain. At initial time, $t = 0$ a jump along $y$ is prescribed to the saturation, then the fluid is allowed to flow forced by the gravity. At $t = 0.3$ we see the impact of the high permeable fracture, that drains the phase 0 from the top and releases it on the bottom tip around $y = 0.3$. As the stationary condition is approached, $t = 2$, the two phases are separated and they occupy half of the domain each.}
    \label{fig:case_1_vertical_saturation}
\end{figure}

The flow in the high-permeability fracture is faster than that in the surrounding materials, so
the fracture drains the heavy fluid from the top and releases it at the bottom. Indeed, observing
Fig.~\ref{fig:case_1_vertical_saturation}, showing the saturation along a vertical line near, but not coincident with, the fracture in the 2D domain, we can see, at time $t = 0.3$, a formation of a local minimum of the saturation located at the top (absorption) and a local maximum at the inner tip of the fracture, at $y=0.3$ (release).

Results regarding the iterative methods, Fig.~\ref{fig:case_1_newton}, are similar to the ones of the previous case.
\begin{table}[h!]
    \centering
    \begin{tabular}{ll}
        \hline
        Matrix intrinsic permeability   & $K_2 = 1$              \\
        Fracture intrinsic permeability & $K_1 = 10$             \\
        Fracture normal permeability    & $k_{\perp,1} = 0.1$    \\
        Fracture cross-sectional area   & $\varepsilon_1 = 0.01$ \\
        Matrix porosity                 & $\Phi_2 = 0.25$        \\
        Fracture porosity               & $\Phi_1 = 0.25$        \\
        Total simulation time           & $ t_{end} = 5$         \\
        Timestep max                    & $\Delta t_{max} = 0.1$ \\
        $E_A$                           & $0.0625$               \\
        \hline
    \end{tabular}
    \caption{Case 1. Vertical. Parameters used for this problem.}
    \label{tab:case_1_vertical_prop}
\end{table}
%
%

\clearpage

\subsubsection{Case 1.c. Slanted fracture, non-conforming grid}
The parameters and properties defining in this case are summarized in Tab.~\ref{tab:case_1_slanted_prop}.

\begin{figure}
    \centering
    \includegraphics[width=0.3\textwidth]{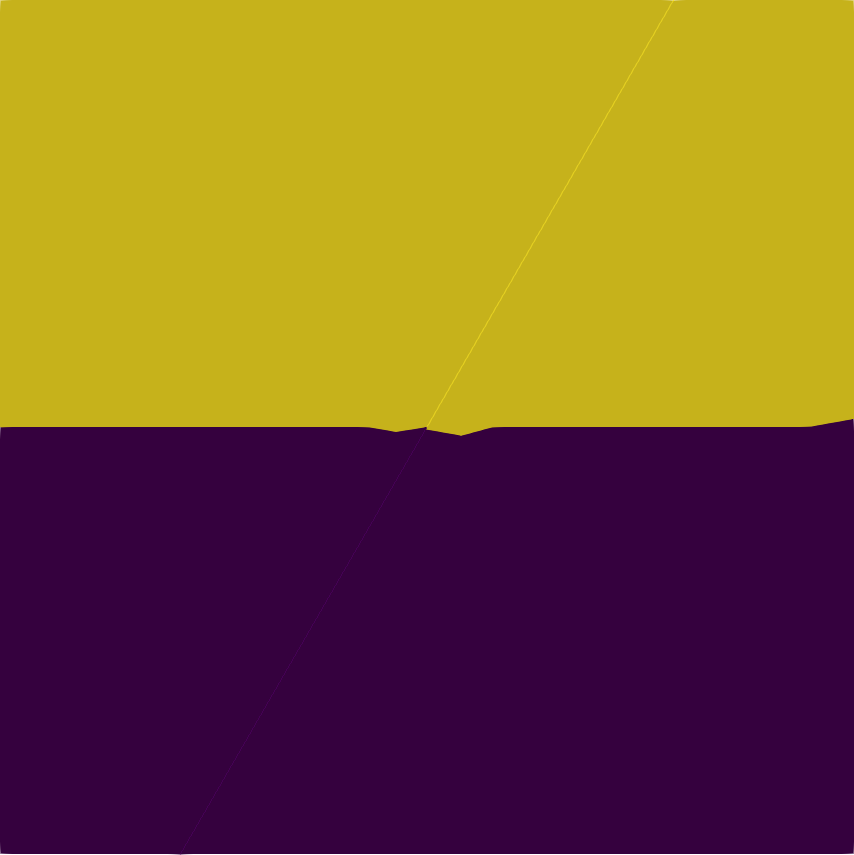}
    \includegraphics[width=0.3\textwidth]{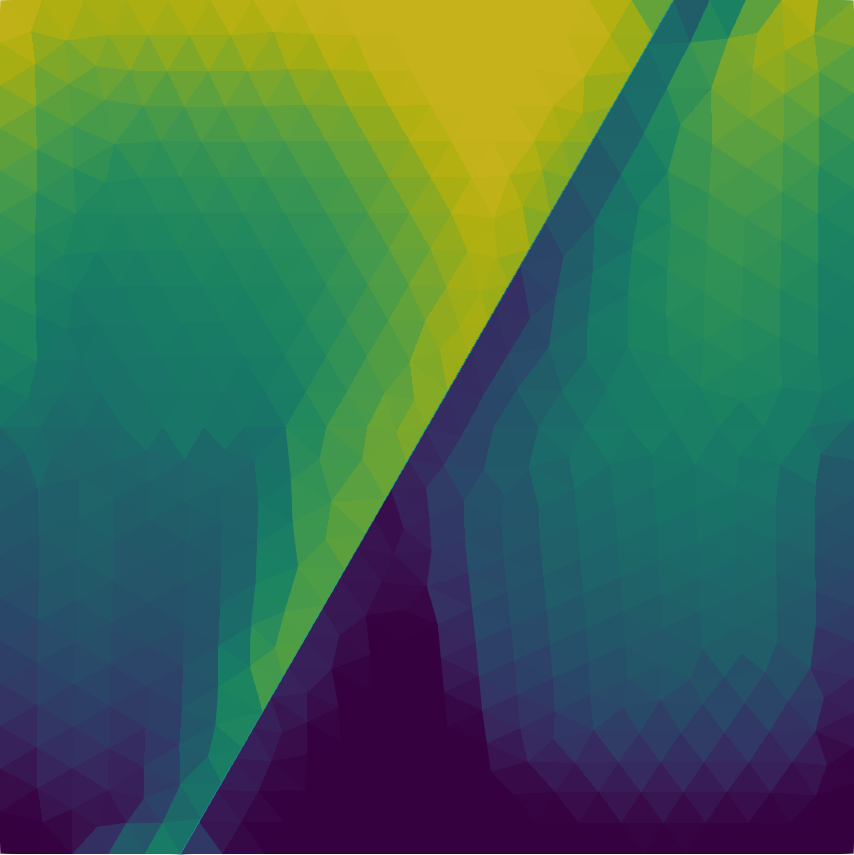}
    \includegraphics[width=0.3\textwidth]{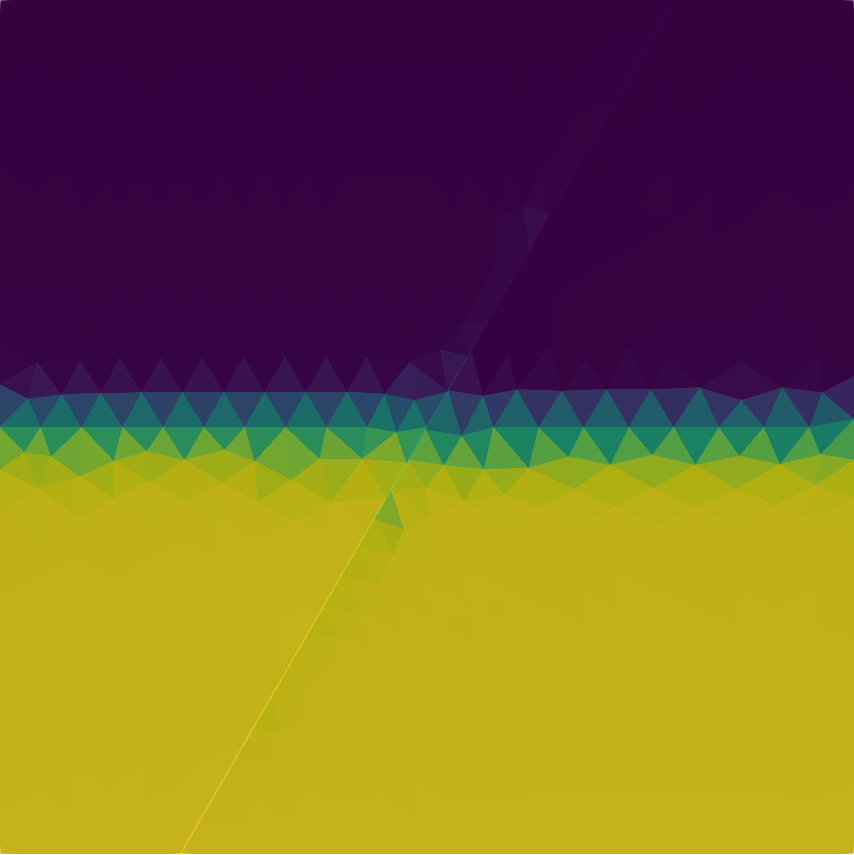} \\
    \vspace{0.1cm}
    \includegraphics[width=0.4\textwidth]{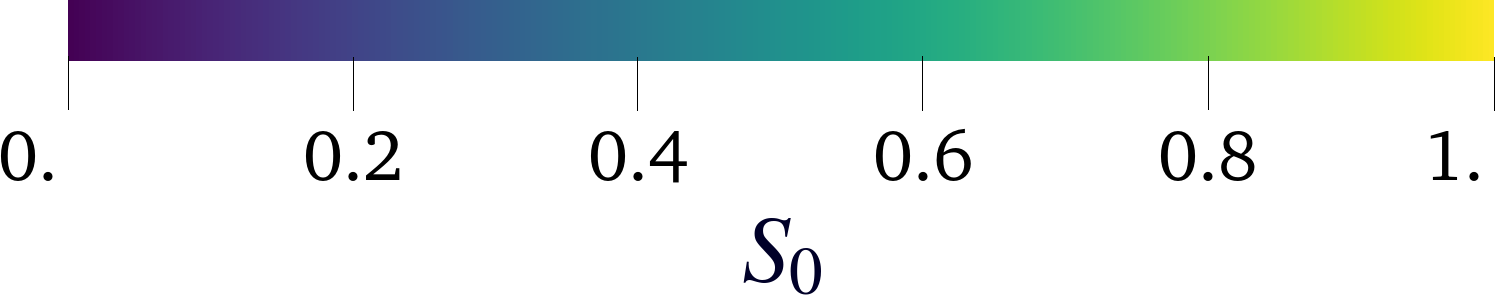}
    \caption{Case 1. Slanted fracture. Three time instances, on the left we have the initial condition in which the heavy phase (denoted with subscript $0$) occupies the upper half of the domain and the light phase occupies the lower part. A non perfectly straight interface is due to the irregularity of the mesh. The middle panel shows the saturation distribution at a later time. We can clearly see the effects of the low-permeable fracture that generates a jump in the saturation value. On the right is the stationary condition, where the phases have swapped positions.}
    \label{fig:case_1_slanted_dyn}
\end{figure}

The normal permeability $k_\perp$ of the fracture is low so it creates an obstacle to the flow through it. Indeed, in Fig.~\ref{fig:case_1_slanted_dyn} we can see that the phases tend to slide along the fracture before reaching the stationary condition where the two phases have swapped position with respect to the initial condition.

\begin{figure}
    \centering
    \includegraphics[width=0.8\textwidth]{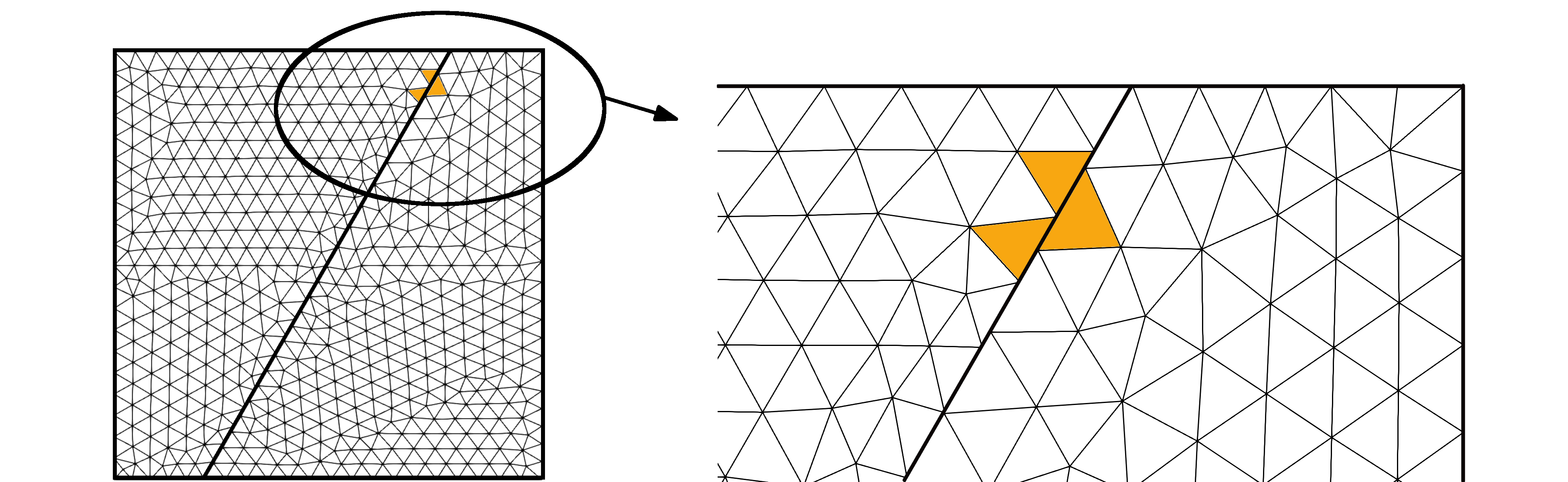}
    \caption{Case 1. Slanted fracture. Zoom on the simplex grid of the 2D domain. The orange elements highlight the lack of conformity at the fracture interface.}
    \label{fig:case_1_slanted_grid}
\end{figure}

We discretized the domain with a simplex non-conforming grid a the interface, as visible in Fig.~\ref{fig:case_1_slanted_grid}. The grid is derived by deforming a conforming grid using a technique that relies on radial basis functions, detailed in \cite{Ballini2024}.
This setup is numerically more challenging the the previous, so a greater gap in the performance is visible in Fig.~\ref{fig:case_1_newton}. For this case, we show also the number of wasted flips of the upwind direction for each subdomain and mortar, bottom panel of Fig.~\ref{fig:case_1_newton}. The 2D and 1D domains exhibit similar trends, with the HU method showing fewer changes in direction. The behavior of the upwind across domains, which is the same for both PPU and HU, varies according to the upwind discretization method used in the neighboring domains and there is not a evident improvement in the HU case. This suggests that a hybrid upwind between domains could lead to advantages.
By comparing these graphs with those showing the cumulative number of Newton iterations, we notice that the Newton iterations are predominantly influenced by the flips in the 2D domain.

Furthermore, increasing the permeability contrast in the 1D and 2D domains does not significantly affect the characteristics of the iterative method. The results are shown in Fig.~\ref{fig:case_1_newton} \textit{Case 1.c} and \textit{Case 1.c hc} (high contrast) where the permeabilities of the fracture are $k_1 = 10^{-6}$ and $K_\perp=10^{-8}$. For \textit{Case 1.c}, the bottom panel of Fig.~\ref{fig:case_1_newton} displays the cumulative number of upwind direction flips for both 2D and 1D domains, as well as for the mortar. 

\begin{figure}[h]
    \centering
    \includegraphics[width=0.32\textwidth]{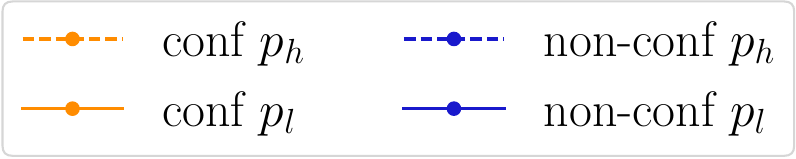}
    \includegraphics[width=0.32\textwidth]{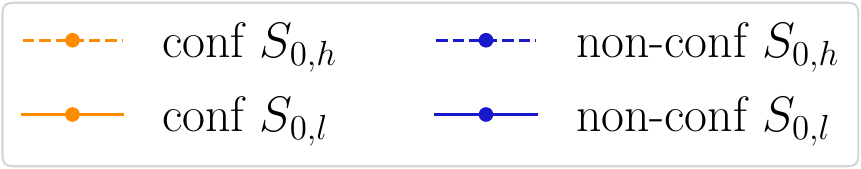}
    \includegraphics[width=0.32\textwidth]{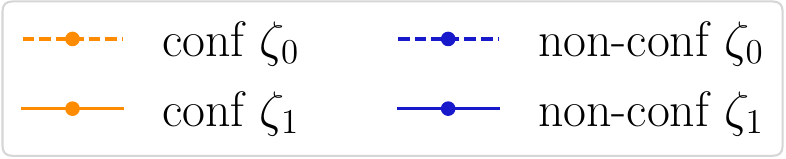}

    \includegraphics[width=0.32\textwidth]{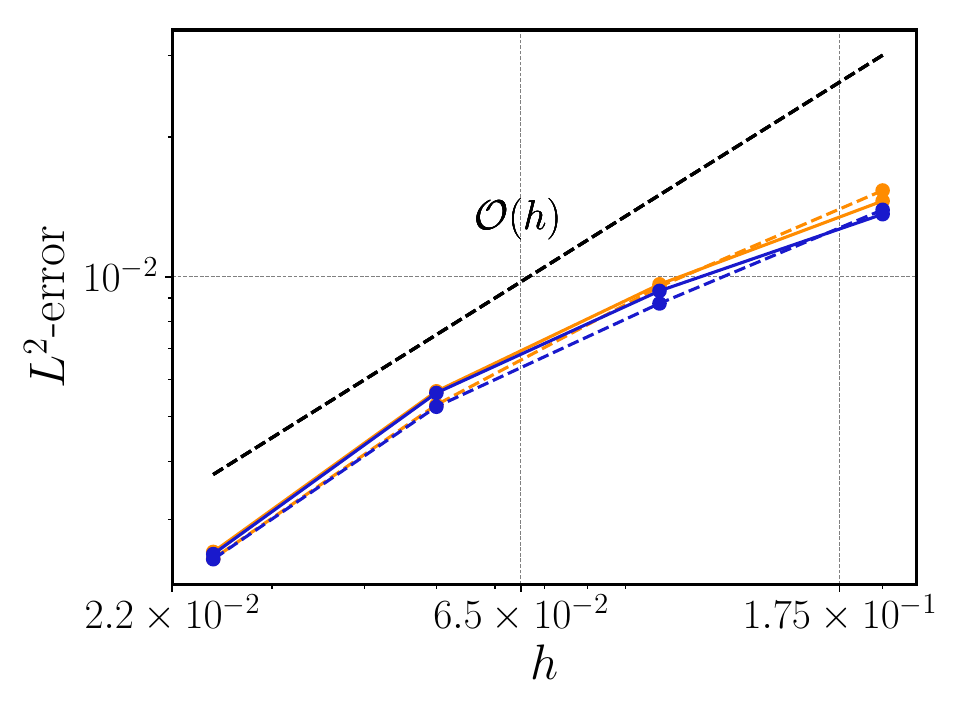}
    \includegraphics[width=0.32\textwidth]{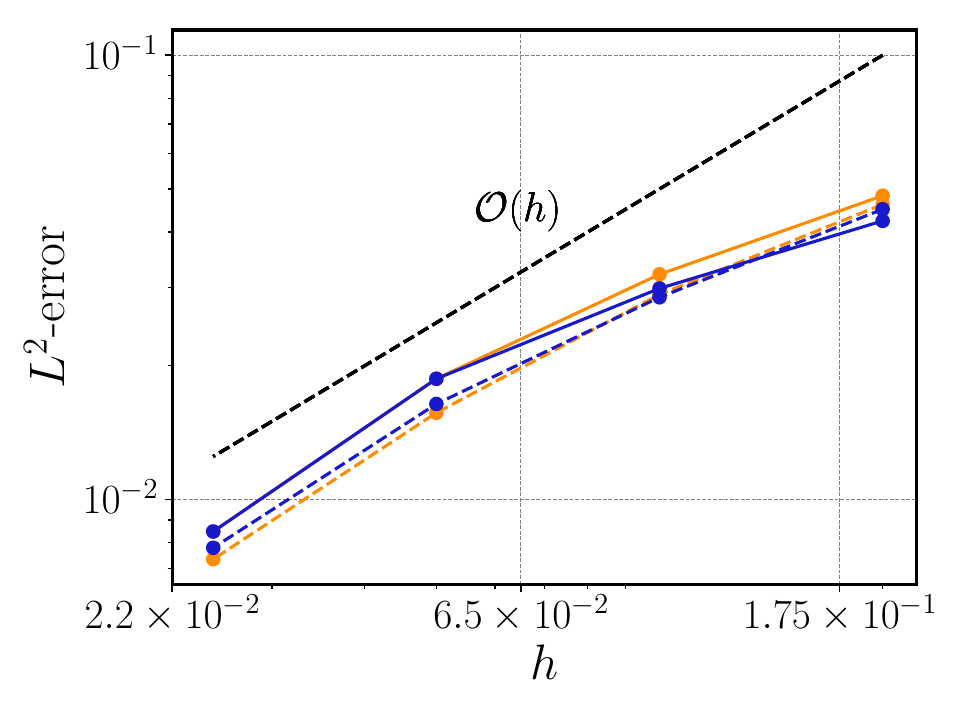}
    \includegraphics[width=0.32\textwidth]{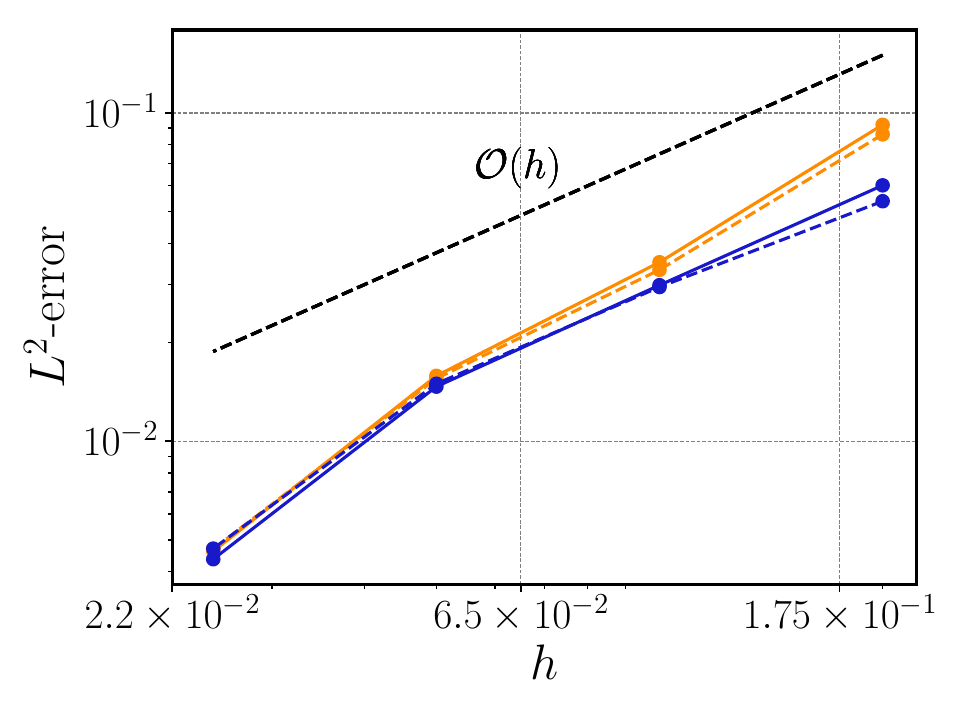}
    \caption{Case 1. Slanted. Spatial convergence on a conforming grid and non-conforming gird.}
    \label{fig:case_1_slanted_convergence}
\end{figure}

We perform a spatial convergence analysis of the HU method, comparing the solution obtained with a reference conforming grid, Fig.~\ref{fig:case_1_slanted_convergence}. The analysis is performed by calculating a single time step from a smooth solution in the primary variables and then calculating the $L^2$-norm of the spatial error with respect to a reference solution computed with the HU on a very fine grid.

The nonconforming grid allows to exploits the full capabilities of the dual-mortar formulation, demonstrating that there are no significant differences on the error using a conforming grid, thus with identity maps, $\Pi$ and $\Xi$, or a nonconforming grid.

\begin{figure}
    \centering
    \includegraphics[width=0.4\textwidth]{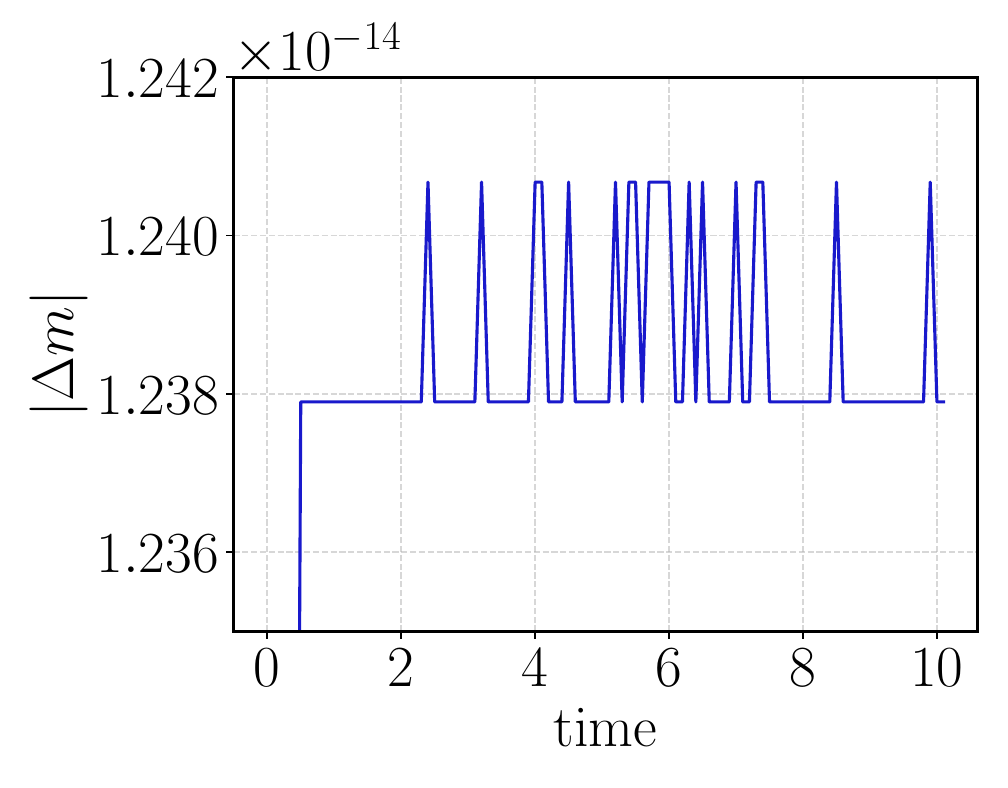}
    \caption{Case 1. Slanted. Variation of the total mass of phase 0 versus time. The variation is referred to the initial mass, so at initial time the variation is $0$ and it does not fit the axis range. The initial mass is $0.127$. Thanks to the finite volume method, the mass is conserved at each timestep.}
    \label{fig:case_1_slanted_mass}
\end{figure}

Thanks to the finite volume method, conservative by construction, the mass of the two phases remains perfectly conserved throughout the entire simulation, as illustrated in Fig.~\ref{fig:case_1_slanted_mass}. Equivalent results are obtained for all the other tests, so they will not be shown.

\begin{table}[h!]
    \centering
    \begin{tabular}{ll}
        \hline
        Matrix intrinsic permeability   & $K_2 = 1$              \\
        Fracture intrinsic permeability & $K_1 = 1$              \\
        Fracture normal permeability    & $k_{\perp,1} = 0.01$   \\
        Fracture cross-sectional area   & $\varepsilon_1 = 0.01$ \\
        Matrix porosity                 & $\Phi_2 = 0.25$        \\
        Fracture porosity               & $\Phi_1 = 0.25$        \\
        Total simulation time           & $ t_{end} = 10$        \\
        Timestep max                    & $\Delta t_{max} = 0.1$ \\
        $E_A$                           & $0.0625$               \\
        \hline
    \end{tabular}
    \caption{Case 1. Slanted fracture. Parameters used in this problem.}
    \label{tab:case_1_slanted_prop}
\end{table}

%


%
%


\subsection{Case 2. Complex fracture network}\label{sec:case_2}
\begin{figure}[h]
    \centering
    \subfloat[Fracture network]{
        \begin{minipage}[t]{0.5\textwidth}
        \centering
        \includegraphics[width=0.5\textwidth]{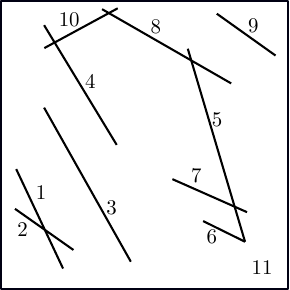}
        \vspace{0.6cm}
        \end{minipage}
    }
    \subfloat[Saturation at $t = 0.013$]{
        \begin{minipage}[t]{0.5\textwidth}
        \centering
        \includegraphics[width=0.5\textwidth]{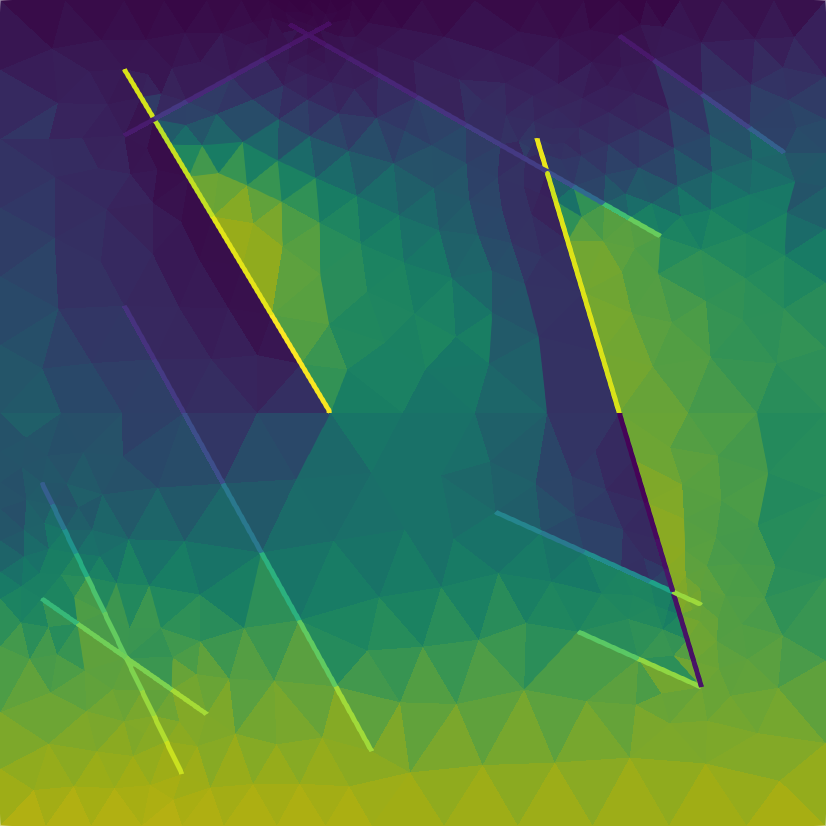} \\
        \includegraphics[width=0.6\textwidth]{case_1_slanted_label.png}
        \end{minipage}
    }
    \caption{Case 2. (a) Domain. The fracture network geometry is taken from \cite{Flemisch2018}.
    (b) Saturation at time $t = 0.013$. The profiles of the saturation in the high-permeable fractures adjust to the surrounding domain. Instead, the low-permeable fracture permit a small motion of fluid across and along them.}
    \label{fig:case_2_domain_and_more}
\end{figure}

The geometry of this test case is taken from \cite{Flemisch2018}.
The presence of many fractures, with intersection of X-type or L-type (fracture 5 and 6) as shown in Fig.~\ref{fig:case_2_domain_and_more}, and the resulting formation of 0D domains, adds a significant challenge to this test. Fractures number 4 and 5 have low permeability, while the other fractures are highly permeable. The permeability of the intersections are set equal to the harmonic average of the intersecting fractures. The parameters and properties defining in this case are summarized in Tab.~\ref{tab:case_2_prop}. 

The low-permeable fractures in addition to creating a barrier, retain the fluid by allowing a slow displacement inside them, as depicted in Fig.~\ref{fig:case_2_domain_and_more}. While in the permeable fractures the saturation values is similar to the one of the surrounding matrix, in the impermeable fractures the saturation distribution is similar to the distribution at the initial time.

\begin{figure}[h]
    \centering
    \subfloat[Cumulative number of flips of the upwind directions.]{
        \begin{minipage}[b]{0.33\textwidth}
        \centering
        \includegraphics[width=1\textwidth]{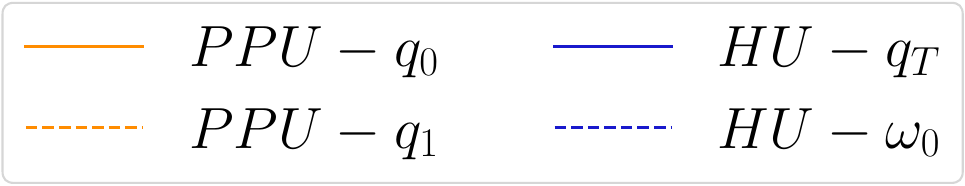}\\
        \includegraphics[scale=0.25]{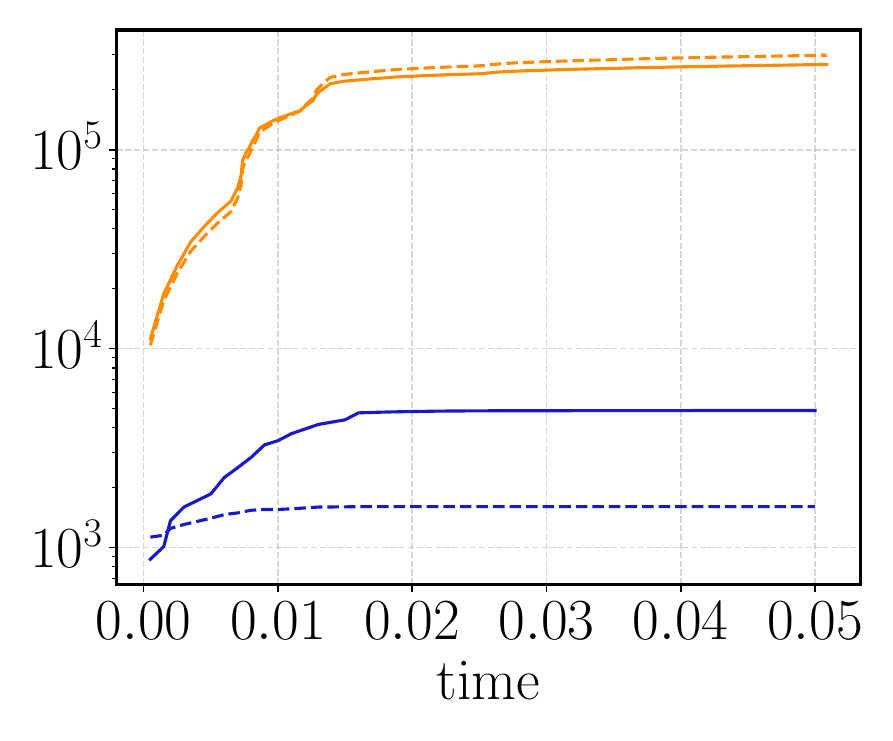} 
        \end{minipage}
    }
    \subfloat[Cumulative number of time cuts.]{
        \begin{minipage}[b]{0.33\textwidth}
        \centering
        \includegraphics[width=0.75\textwidth]{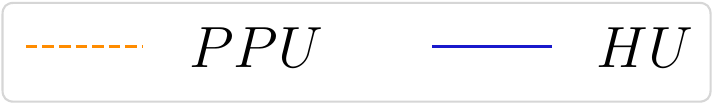}\\
        \includegraphics[scale=0.25]{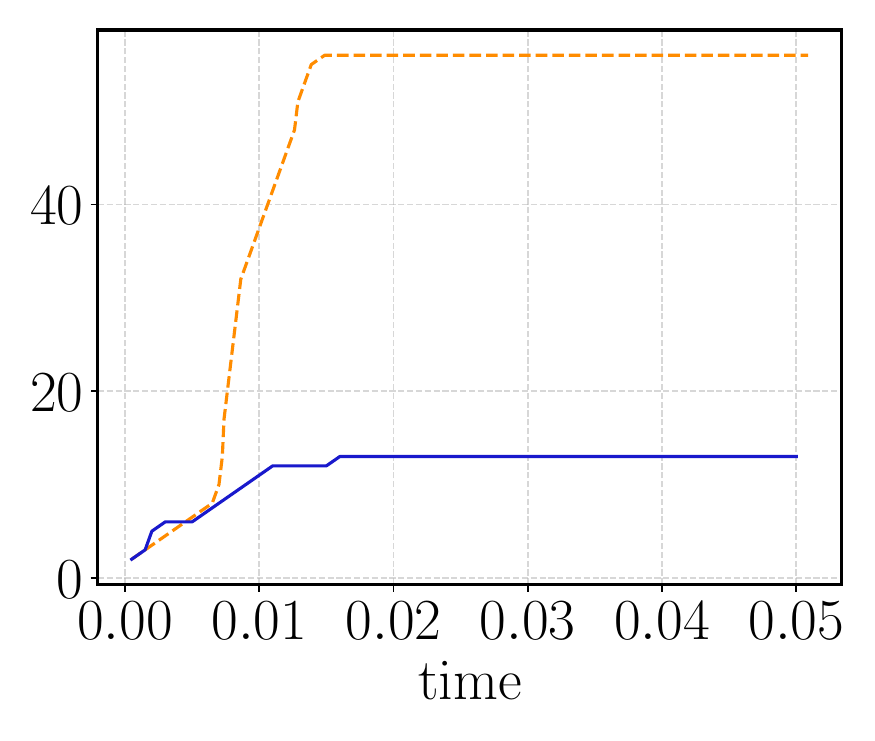} 
        \end{minipage}    
    }
    \subfloat[Cumulative number of Newton iterations.]
    {
        \begin{minipage}[b]{0.33\textwidth}
        \centering
        \includegraphics[width=0.75\textwidth]{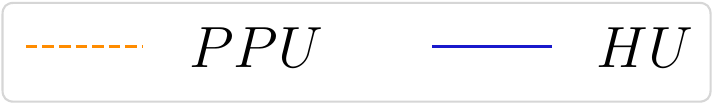}\\
        \includegraphics[scale=0.25]{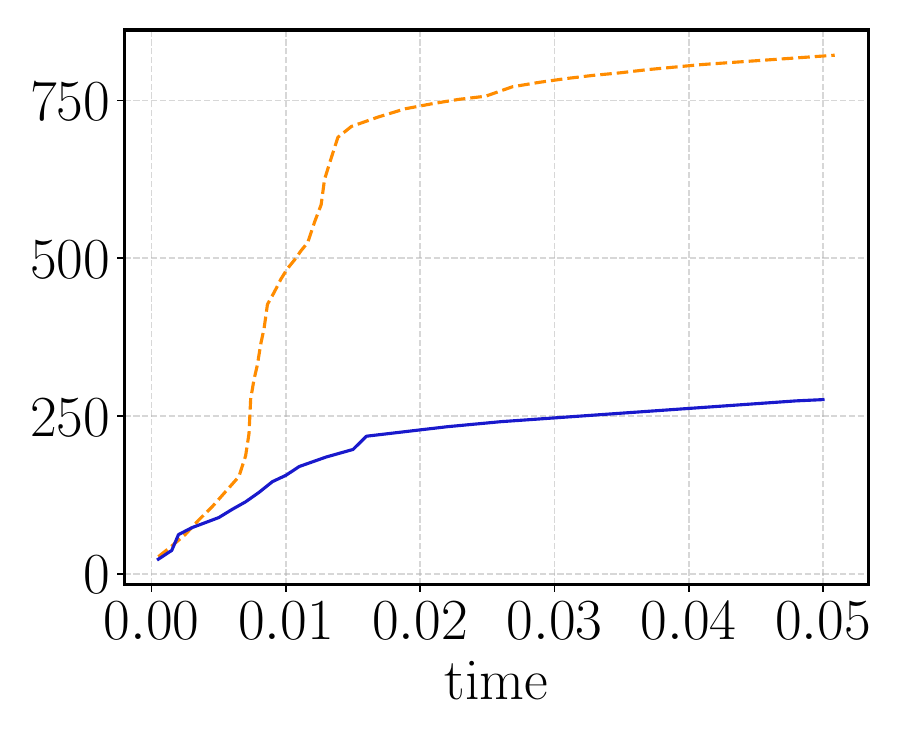} 
        \end{minipage}    
    }
    \caption{Case 2. Characteristics regarding the iterative method.}
    \label{fig:case_2_newton}
\end{figure}

The cumulative number of flips of the upwind direction, the number of timestep reductions  and cumulative Newton iteration, Fig.~\ref{fig:case_2_newton}, show the better performance of HU compared to PPU, in particular, the latter requires a number of timesteps cuts up to 4 per time step, causing a large amount of wasted iteration. Indeed, the end of the simulation, the number of iterations required by HU is approximately three times lower. The simulation thus showed that HU significantly outperforms PPU also for more complex geometries.

\begin{table}[h!]
\centering
\begin{tabular}{ll}
\hline
    Matrix intrinsic permeability & $K_{11} = 100$ \\
    Fracture $i \neq 4,5$ intrinsic permeability & $K_i = 100$ \\
    Fracture 4 intrinsic permeability & $K_4 = 0.01$ \\
    Fracture 5 intrinsic permeability & $K_5 = 0.01$ \\
    Fracture $i \neq 4,5$ normal permeability & $k_{\perp,i} = 100$ \\
    Fracture 4 normal permeability & $k_{\perp,4} = 0.01$ \\
    Fracture 5 normal permeability & $k_{\perp,5} = 0.01$ \\
    Intersection normal permeability & $k_{\perp,i} = \frac{1}{1/k_{\perp,j} + 1/k_{\perp,q}}$ \\
    Fracture cross-sectional area & $\varepsilon_i = 0.01$ \\
    Intersection cross-sectional area & $\varepsilon_i = 0.01$ \\
    Matrix porosity & $\Phi_{11} = 0.25$ \\
    Fracture porosity & $\Phi_i = 0.25$ \\
    Intersection porosity & $\Phi_i = 0.25$ \\
    Total simulation time & $ t_{end} = 0.05$ \\
    Timestep max & $\Delta t_{max} = 2\times 10^{-3}$ \\
    $E_A$ & $6.25$ \\
\hline
\end{tabular}
\caption{Case 2. Parameters used for this problem. The intersection normal permeability is the harmonic average of the intersecting fractures $\Gamma_j$ and $\Gamma_q$.}
\label{tab:case_2_prop}
\end{table}

%


\subsection{Case 3. Network with small features}\label{sec:case_3}
\begin{figure}[h]
    \centering
    \subfloat[Fracture network]{
    \includegraphics[width=0.22\textwidth]{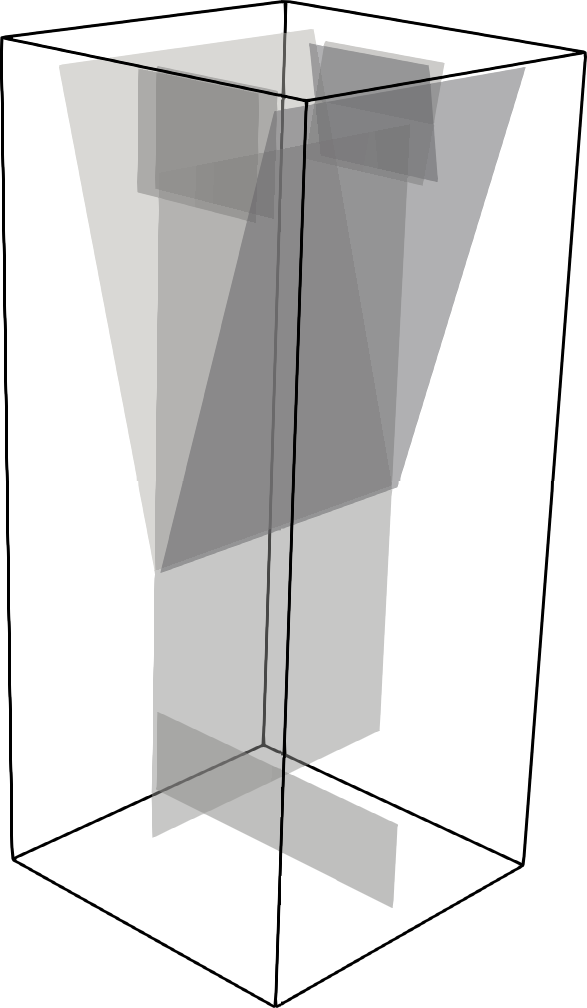}
    }
    \subfloat[Ill-shaped cells]{
    \includegraphics[width=0.26\textwidth]{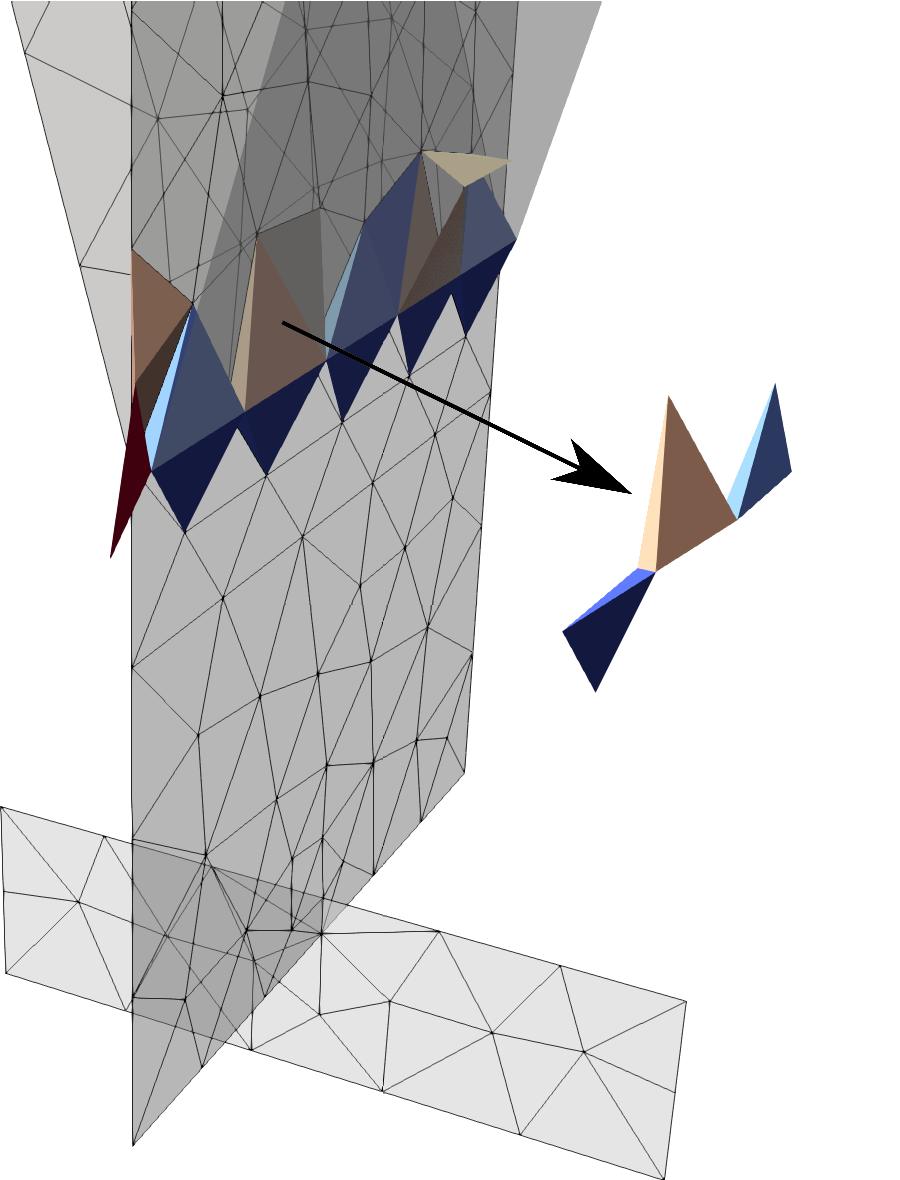}
    }
    \subfloat[$t=1.3\times10^{-3}$]{
    \includegraphics[width=0.22\textwidth]{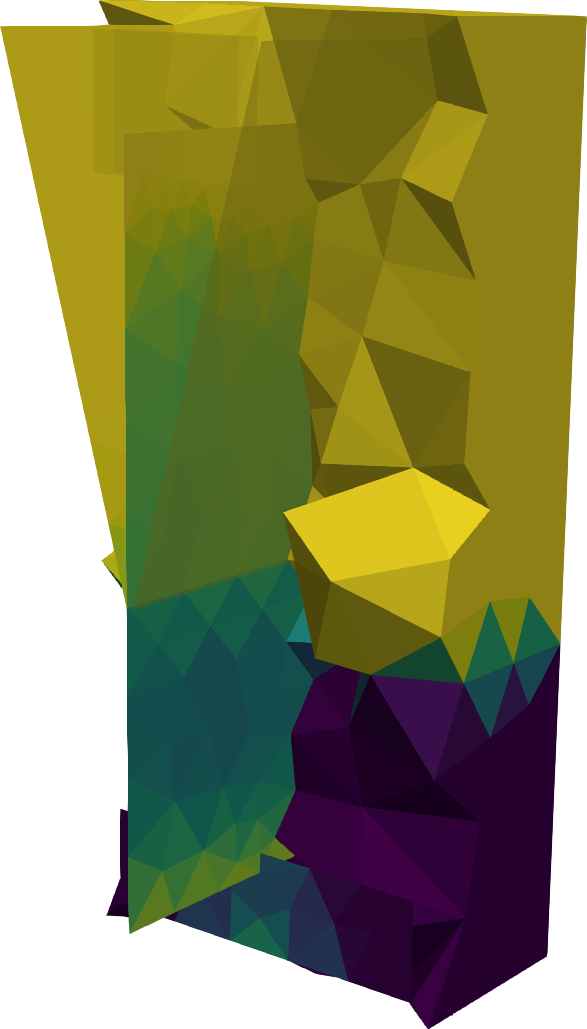}
    }
    \subfloat[$t=2.86\times10^{-3}$]{
    \includegraphics[width=0.22\textwidth]{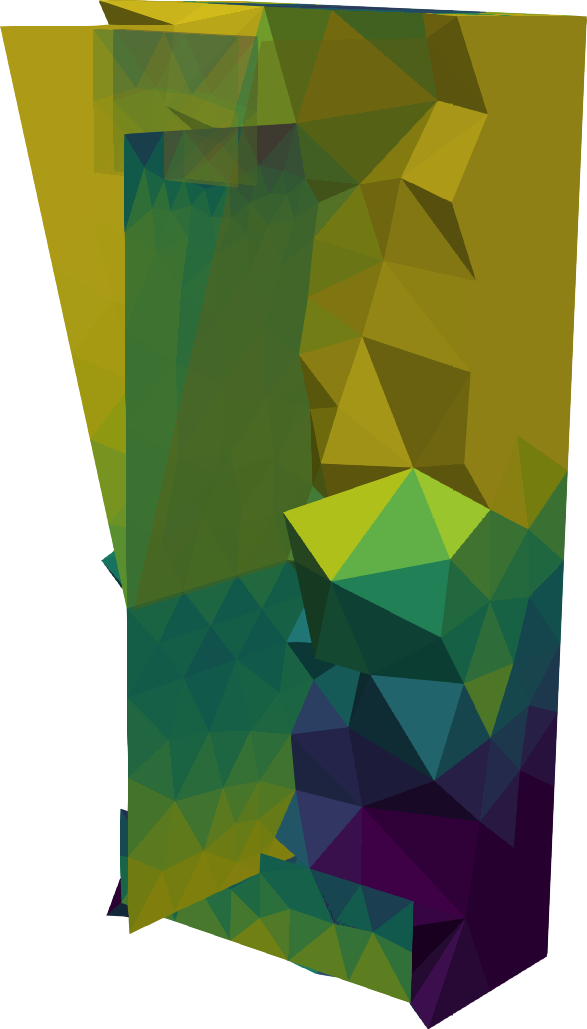}
    }
    \caption{Case 3. (a) Domain and fracture network. 
    (b) The fracture configuration with narrow angles leads to ill-shaped grid elements. The highlighted elements have a ratio of the circumscribed and inscribed sphere radius between $3.5$ and $4.1$.
    (c) and (d) Saturation at time $t = 1.3\times10^{-4}$ and $t = 2.86\times10^{-3}$, respectively. We can appreciate different time scales due to different $E_A$ numbers, specifically $11.25$ in the matrix and $875$ in the fracture, the fast one inside the fractures and the slow one in the 3D domain. Some fractures are made partially transparent for graphical reason.} 
    \label{fig:case_3_domain_and_more}
\end{figure}

The geometry and the domain are replicated from \cite{Berre2021}. We list in Tab.~\ref{tab:case_3_prop} the main properties. For this case, we set the Newton tolerance $tol = 2\times10^{-5}$. The fracture network is high permeable whereas the surrounding matrix is low permeable. Fig.~\ref{fig:case_3_domain_and_more} illustrates the 2D fracture network immersed in a 3D domain. The fracture intersection are of X and Y-type with small angles. This complicated geometry entails the generation of ill-shaped grid elements, a challenge for the discretization methods, Fig.~\ref{fig:case_3_domain_and_more}. 

As in the previous case, we set the heavy phase to lay initially in the top part of the domain and the motion is then forced by the gravity. The permeability contrast produces a fast dynamics in the fracture that lasts till around time 0.1 before a slow motion takes place in the whole domain till the end of the simulation, we can appreciate the different speed of diffusion inside the fracture and in the matrix due to different $E_A$ numbers, $11.25$ for the matrix and $875$ for the fracture. In Fig.~\ref{fig:case_3_domain_and_more} the saturation during the transient is shown. A small diffusion is visible in the 3D domain, while a large displacement of the phases occurred in the fracture network.
Results in Fig.~\ref{fig:case_2_newton} show that the PPU fails to converge at the beginning of the simulation, during the fast dynamics, while the increased robustness of HU enabled this challenging simulation to be completed with a limited number of time step cuts.

\begin{table}[h!]
\centering
\begin{tabular}{ll}
\hline
    Matrix intrinsic permeability & $K_3 = 100$ \\
    Fracture $i$ intrinsic permeability & $K_i = 10^{4}$ \\
    Intersection $i$ intrinsic permeability & $K_i = 10^4$ \\
    Fracture $i$ normal permeability & $k_{\perp,i} = 10^{4}$ \\
    Intersection $i$ normal permeability & $k_{\perp,i} = 10^4$ \\
    Fracture cross-sectional area & $\varepsilon_i = 0.01$ \\
    Intersection cross-sectional area & $\varepsilon_i = 0.01$ \\
    Matrix porosity & $\Phi_{3} = 0.2$ \\
    Fracture porosity & $\Phi_i = 0.2$ \\
    Intersection porosity & $\Phi_i = 0.2$ \\
    Total simulation time & $ t_{end} = 0.01$ \\ 
    Timestep max & $\Delta t_{max} = 10^{-5}$ \\
    Newton tolerance & $tol = 2\times10^{-5}$ \\ 
    $E_A$ & $11.25$ \\
\hline
\end{tabular}
\caption{Case 3. Parameters defining the test case.}
\label{tab:case_3_prop}
\end{table}
\begin{figure}[h]
    \centering
    \subfloat[Cumulative number of flips of the upwind directions.]{
        \begin{minipage}[b]{0.33\textwidth}
        \centering
        \includegraphics[width=1\textwidth]{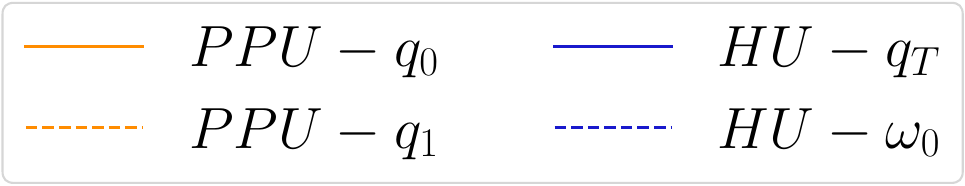}\\
        \includegraphics[scale=0.25]{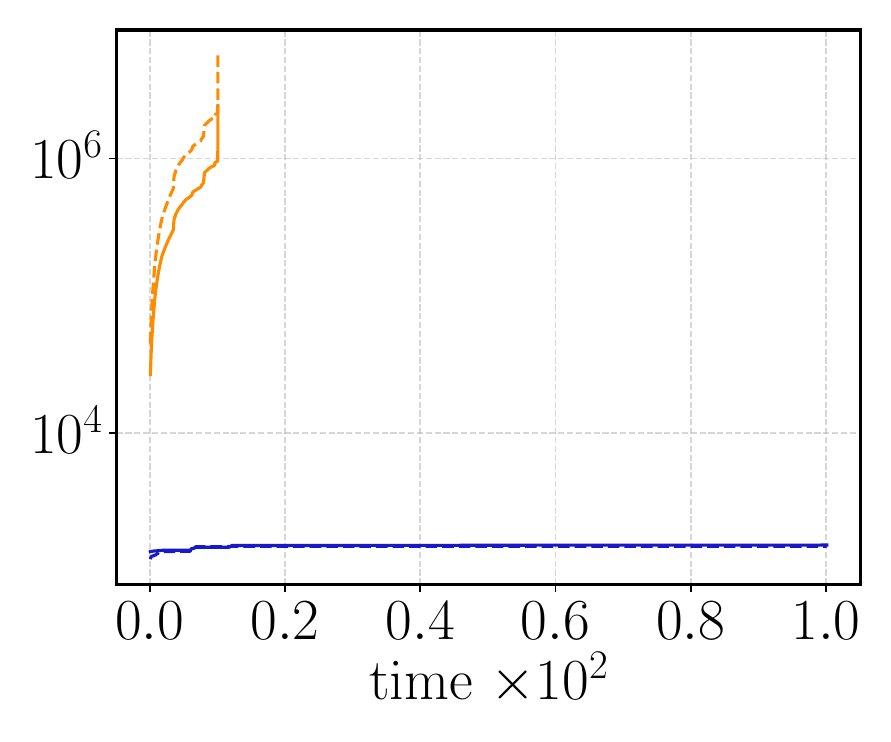} 
        \end{minipage}
    }
    \subfloat[Cumulative number of time cuts.]{
        \begin{minipage}[b]{0.33\textwidth}
        \centering
        \includegraphics[width=0.75\textwidth]{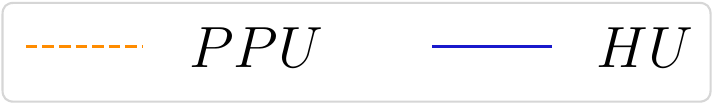}\\
        \includegraphics[scale=0.25]{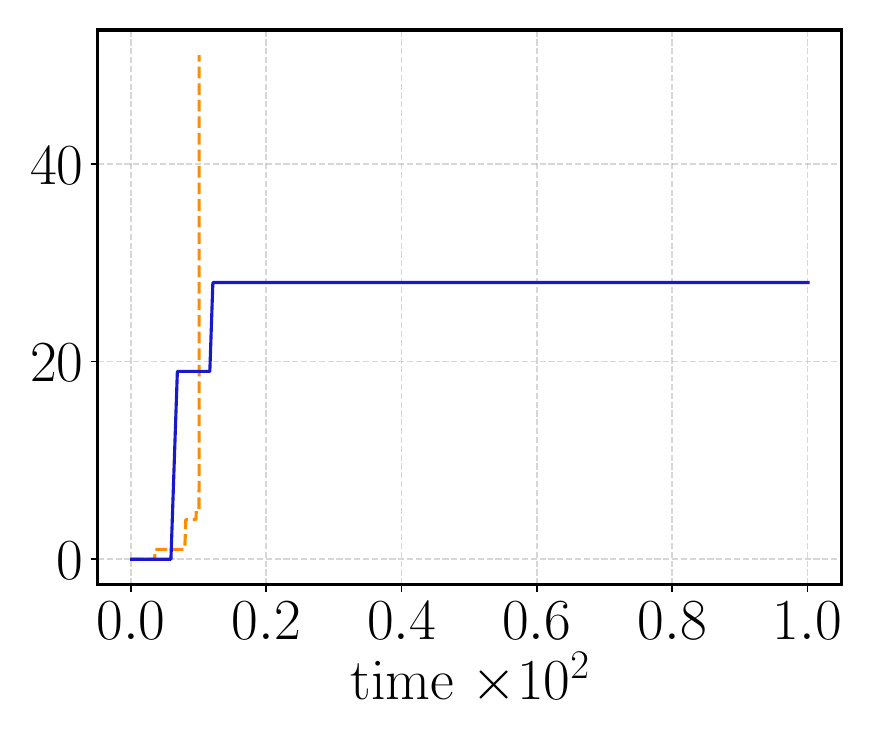} 
        \end{minipage}    
    }
    \subfloat[Cumulative number of Newton iterations.]
    {
        \begin{minipage}[b]{0.33\textwidth}
        \centering
        \includegraphics[width=0.75\textwidth]{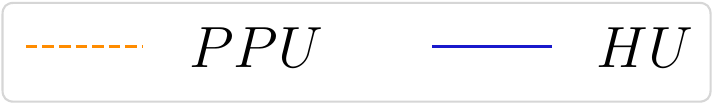}\\
        \includegraphics[scale=0.25]{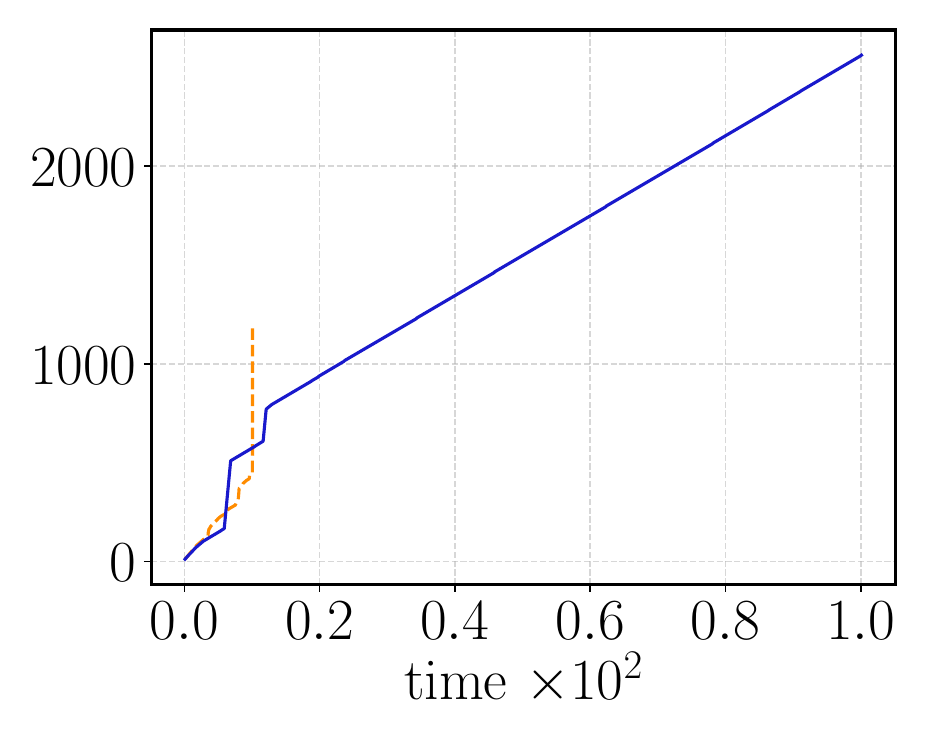} 
        \end{minipage}    
    }
    \caption{Case 3. Characteristics regarding the iterative method. The PPU method fails to converge up to a reasonable timestep size while the HU proves to be more robust.}
    \label{fig:case_3_newton}
\end{figure}
%


\clearpage

\section{Conclusion}
\label{sec:conclusion}
We address the physical problem of two-phase flow in the subsurface with a fractured matrix rock. We model fractures, and possibly intersections, through a dimensionally reduced object obtaining a mixed-dimensional domain. The fluid motion is governed by partial derivative equations defined in each subdomain, while the interaction between the subdomains is described by Lagrange multipliers that are represented by fluid fluxes. The fluid flow is dictated by possible sources/wells, boundary conditions and, most importantly, gravity. This configuration is known to generate countercurrent flows, a condition that spoils the performance of the nonlinear solver (Newton).

We tackle the convergence problem by working on the discretization of the fluxes. We then extended the work done in \cite{Bosma2022} on a hybrid upwind strategy to the case of a mixed-dimensional framework. We implement the method in PorePy, a simulation tool for fractured and deformable porous media suitable for the mixed-dimensional problem \cite{Keilegavlen2019}.

We test the discretization method on three different geometries, both 2D and 3D, with intersecting fracture network, discretized with simplex and hexahedral meshes, and we test the method on different flow regimes, imposed by different rock properties. 
We show numerically the convergence of the method, even with non-conforming grids at subdomain interfaces. In each test case, the proposed method reduces the number of Newton iterations. In particular, in the third case (\Cref{sec:case_3}), the standard discretization method fails to converge, making it impossible to complete the time-dependent simulation. On the other hand, an increase in numerical diffusion is observed.

Given the promising results obtained, further developments will be undertaken, such as, from the physical point of view, extension to $n$-phases, high compressible fluids, inclusion of capillary effects and chemical reactions. From the numerical point of view, a hybrid upwind strategy should also be adopted at the interfaces to further improve performance, in addition to a study aimed at decreasing the numerical diffusion added by the discretization method. Moreover, despite an increment in the computational cost for solving the linear system, we expect an improvement in the accuracy by using of MPFA and gravity consistent transmissibilities in \eqref{eq:q_ell} instead of TPFA.

The numerical scheme has been shown to be robust and effective in reducing the number of Newton iterations, resulting in benefits in computational cost of the simulation. 

\section*{Acknowledgement}
The authors declare that they have no known competing financial interests or personal relationships that could have appeared to influence the work reported in this paper. 

The present research is part of the activities of ``Dipartimento di Eccellenza 2023-2027'', Italian Minister of University and Research (MUR), grant Dipartimento di Eccellenza 2023-2027, under the project ``PON Ricerca e Innovazione 2014-2020".

The authors gratefully thank Hamon Fran\c{c}ois for the insightful discussions.

\small
\bibliography{all}
\bibliographystyle{abbrv}

\end{document}